\documentclass[final]{siamltex}

\usepackage{srcltx}

\usepackage{graphics}
\usepackage{epsfig}
\usepackage[mathscr]{eucal}
\usepackage{amscd}
\usepackage{amsfonts}
\usepackage{amsmath,  amssymb}
\usepackage{latexsym}
\usepackage[usenames,dvipsnames]{color}
\usepackage{amsmath}
\newcommand{\bigO}[1]{\mathcal{O}(\varepsilon^#1)}
\newcommand\coll{u}
\newcommand\normal{n}
\newcommand\dep{v}
\newcommand\Diff{d}
\newcommand\dca{a}
\newcommand\dcb{b}

\usepackage{xargs}
\usepackage{etoolbox}
\newcommandx{\inRange}[2][2=]{\ifstrequal{#2}{}
  {\in\{1,\ldots,#1\}}
  {\in\{#1,\ldots,#2\}}
}
\def\vec{\overrightarrow}
\def\d{\delta}

\def\be{\begin{equation}}
\def\ee{\end{equation}}
\def\bse{\begin{subequations}}
\def\ese{\end{subequations}}
\def\bge{\begin{eqnarray}}
\def\bgee{\begin{eqnarray*}}
\def\ege{\end{eqnarray}}
\def\egee{\end{eqnarray*}}

\usepackage[pagewise]{lineno}
\title{Colloidal Transport in Locally Periodic Evolving Porous Media -- An Upscaling Exercise}

\author{ Adrian Muntean 
         \thanks{ Department of Mathematics and Computer Science, Karlstad University, Sweden ({\tt adrian.muntean@kau.se}). }
 \and
Christos V. Nikolopoulos       
         \thanks{Department of Mathematics,
University of the Aegean, Gr-83200 Karlovassi, Samos, Greece. 
({\tt cnikolo@aegean.gr}).} 
}       

\date{}
\begin{document}
%_________________________________________________________________
\maketitle

\begin{abstract}
We derive an upscaled  model describing the aggregation and deposition of colloidal particles within a porous medium allowing for the possibility of local clogging of the pores. At the level of the pore scale, we extend an existing model for colloidal dynamics  including the  evolution of free interfaces separating  colloidal particles deposited on solid boundaries (solid spheres) from the colloidal particles transported through the gaseous parts of the porous medium.  As a result of deposition, the solid spheres grow reducing therefore the space available for transport in the gaseous phase.  
 Upscaling procedures are applied and several classes of macroscopic models together with effective transport tensors are obtained, incorporating explicitly the local growth of the solid spheres. The resulting models are solved numerically and various simulations are presented. In particular, they are able to detect clogging regions, and therefore, can provide estimates on the storage capacity of the porous matrix.
%n unsaturated
\end{abstract}

\vspace{0.5in}
\begin{keywords} 
Evolving porous media, \, Colloidal dynamics, \, Clogging,    Storage capacity, \, Asymptotic homogenization,  Numerical simulations.
\end{keywords}
%
%\vspace{0.5in}
\begin{AMS}
Primary 35B27, 76M50; Secondary 76R50, 76VXX, 76M10 %keywords checked AM
\end{AMS}

%\vspace{0.5in}

%________________________________________________________________
\section{Introduction}

We investigate the effect the motion, aggregation, fragmentation and deposition of colloidal particles through porous materials  have on effective transport properties. Colloids are very small particles (from 1 to 100 nm) that play nowadays a very important role in the success of biological and technological applications, like waste water treatment, food industry, 3D printing, drug-delivery design, etc. A mathematically interesting aspect which complicates also the physics of the situation is that  colloids like to aggregate and build large clusters which then fragmentate if some critical thresholds are met. We refer the reader to, for instance, 
 \cite{Elimelech, Peukert, Camejo},   and \cite{Bart} for more details on the mathematical modelling of the aggregation and coagulation processes.   

The main subject of this work is to extend the modelling approach presented in \cite{Krehel1,  Krehel2} (see also in the references cited therein)
 to include the gradual increase of the solid skeleton of the  porous medium due to colloids deposition via a moving boundary modelling strategy. More specifically, we model a reaction-diffusion system operating within a porous medium
consisting of colloidal aggregates of varying size. These aggregates (populations of clusters) interact with each other and with the porous matrix. 
As main effect of the colloidal aggregation and of the deposition of the large clusters  within the evolving porous network, either the transport, or the storage capacity or both are reduced. To explore the multiscale factors behind this effect, we wish to capture {\em how} the local blocking of the pores (the so-called {\em pore clogging}) affects the effective transport and storage. The clogging effect will be pointed out in terms of closed-form formulas for all  effective coefficients. Obtaining quantitively correct entries in the tortuosity tensor requires a deep understanding on how ink-bottle pores and blind pores become closed pores or on how through pores clog due to a heterogeneous deposition of colloids.  Under some  assumptions, these entries will be obtained in this framework through a locally-periodic homogenization asymptotics as in \cite{Chechkin,fatima2011homogenization,Zhang} combined with a suitable matched-asymptotics approach as in \cite{Lacey2,nikolopoulos2004model}.  The numerical evaluation of these tensors will illustrate the appearance of clogging. The situation remotely resembles the filter blockage case discussed in \cite{Dalwadi-Bruna2}. Using a similar methodology as in \cite{Krehel1}, one can also point out the effect of a nonlinear dynamic blocking function on the colloidal deposition.  

It is worth noting that such investigations of the occurrence of clogging are  relevant in a number of engineering directions, especially as {\em multiscale storage} of unwanted particles (chemical species, defects, mechanical damage, etc. ) is concerned. For instance, the efficient and safe storage of injected supercritical carbon dioxide (${\rm CO_2}$) underground is now one potential solution for reducing ${\rm CO_2}$ emissions in the atmosphere; see e.g.  \cite{Rimmele} or \cite{Nordbotten}. Interface and grain boundaries have shown to be ideal sinks for radiation-induced defects. Especially in nano-structured materials, such interfaces can become ideal places for Helium storage; compare e.g. \cite{Beyerlein,Zhu}. Yet another typical research field fitting well this problem setting is filtration combustion \cite{Toledo}.  As one can expect after weighing the number of mathematical difficulties in handling simultaneously the multiscale (often stochastic) structure of the material, the motion of the many microscopic interfaces, the strong coupling and the involved nonlinearities in the physical and chemical process of the setting, there is not yet a consensus on how to approach the question of averaging  reactive, actively evolving 
 porous media. A couple of attempts have been successful in dealing with particular moving interface scenarios. Concerning the employed mathematical techniques, we refer the reader to averaging approaches involving a direct use of suitable level set equations (cf. e.g. \cite{Tycho,ray2013drug,ray2015}), phase field equations (cf. e.g. \cite{Eck}), direct handling of the moving boundary motion e.g. by a precise accounting of free-interfaces-induced volume changes (cf. e.g. \cite{Malte}), or of local balance laws written in terms of a time and space-dependent porosity seen as  an intensive quantity (cf. e.g.\cite{Bruna-Chapman,Dalwadi-Bruna1,Wilmanski}). 
 
 %________________________

Summarizing,  in the present  work we derive an upscaled  model describing the aggregation and deposition of colloidal particles 
within a porous medium with growing spheres as internal microstructure. Two upscaling procedures are applied and several classes of macroscopic models together with effective transport tensors are obtained. They all depend on the local growth of the solid spheres. The simulations based on our models allow us to capture the phenomenon of clogging and to estimate a posteriori how much mass can be stored in the porous matrix. 
 
%_________________________

The paper is structured as follows: The geometry of the porous material as well as the setting of the equations for the microscopic (pore-level) model are given in Section  \ref{micro-eqs}. Section \ref{nondim} contains the non-dimensionalization of the system, while in Section \ref{LPH} and in \ref{large} we provide the details of the upscaling procedures for a couple of different scenarios allowing for a locally periodic distribution of pores (perforations) as well as for the possibility of touching of the microstructures.  Moreover  a  finite element  scheme for the numerical approximation of the system  is presented in Section \ref{simulation}.
 A couple of numerical simulations illustrate the behavior of the Galerkin approximations of the  solution to the upscaled model. The paper concludes with the discussion from Section \ref{discussion}.

%_________________________________________________________________________________
%__________________________________________________________________________________

\section{The microscopic model}\label{micro-eqs}

The colloidal population consists of identical particles called primary particles. They aggregate to bigger shapes and then they fragmentate and finally reach some size, say $i$.  
We refer to each particle of size $i$ as a member of the $i^{th}$ species (or $i^{th}$ cluster). Each colloidal species is characterized by the number of primary particles that  it contains. We have in view here only those situations where each considered size contains a huge number of particles. To describe the situation, we introduce the notation  
$u_i$ when referring to the population of particles of size $i$, $i\in \{1,\ldots, N\}$ for a fixed choice  $N$. The value $N$ corresponds to the biggest allowed size.

\paragraph{Description of the porous domain}

 The porous matrix we have in mind is assumed to have  a locally periodic structure as described in Figure~\ref{fig:locally-periodic}, which resembles the ordered porous medium from \cite{Bruna-Chapman}.
   A two-dimensional microstructure approach is adapted here corresponding to a pseudo 3D dimensional geometry,
   e.g. see Figure \ref{cellfig} 
    (This could be  reduced to a even simpler and less realistic one-dimensional microstructure approach as in ~\cite{Lacey1}).  
 More precisely,  $\Omega$ is a bounded  domain in $\mathbb{R}^2$ representing our porous material\footnote{A similar discussion can be done for the case of $\Omega\subset\mathbb{R}^3$ with a minimal number of modifications mostly concerning scaling factors. }. 
This is assumed to consist of an infinite number of identical cells, each one corresponding to a single pore of the material.  
% Also $\Omega$ will be the bounded domain of $\mathbb{R}^n$ that represents a single cell.
More specifically,  for $\ell>0$  we take
 $Y=\{ \sum_{i=1}^2 \lambda_i \vec{e}_i \,:\, 0<\lambda_i<\ell\}$ an $\ell$-cell in $\mathbb{R}^2$ representing a single pore of the material,
  $ Y_1(x)\subset Y$ open, representing the solid grain and $ Y_0=Y\setminus clos(Y_1)$ the rest of the cell. The boundary of $ Y_1(x)$ will be denoted by
  $\Gamma(0,x)=\partial Y_1(x)$ and is assumed to be piecewise smooth.
  %___________________________________
%   In addition we have 
%  $\Omega^\varepsilon=\Omega\setminus \Omega_0^\varepsilon $ the matrix skeleton in $\Omega$,  
% $\Omega_0^\varepsilon=\cup\{ \varepsilon Y_0^k \, : \, Y_0^k\subset \Omega^\varepsilon \}$ the array of pores (voids)
% while $\Gamma^\varepsilon=\partial \Omega_0^\varepsilon $ represents the pores boundary.
% The colloidal species can be either mobile, hosted in $\Omega_0^\varepsilon$, or immobile, sedimented on
%$\Gamma^\varepsilon$.
%_____________________________________ 
As the process evolves, the mobile colloidal species which reach a certain cluster size tend to deposit at the boundary $\Gamma(0,x)$ of the solid grain $Y_1(x)$.
The basic picture we have in mind is that the solid grain is a disc of radius $r$. As more species deposit on
its boundary,  its radius gradually increases giving rise to a ball $B(r)$ with $ Y_1=B(0)\subset B(r)$.
 We use a radially symmetric setting by assuming that the accumulation of the colloidal species upon the solid grain is uniform.
 We point this out using the notation $\Gamma=\Gamma(r)$, meaning in fact $\Gamma=\Gamma(r(t,x))$. 
 We also need to emphasize that an additional assumption  used in the above setting is that the shape of the colloidal particles  (spherical) does not play an important role in the accumulation process. Inclusion of this characteristic together also with gravitational effect for large $N$ are not taken into account for simplicity reasons.

\paragraph{Setting of the model equations in the patch $Y\subset \Omega$}

As the mobile matter deposits on $\Gamma(r)$, the corresponding deposited mass takes up
volume and leads to the increase in $r$, resulting in a gradual decrease in the porosity of the material, visible at the
macroscopic scale.  To describe the deposition process we use a linear evolution equation and call $v$ the mass of all deposited colloidal populations. The mobile species deposited on $\Gamma(r)$ 
 transform into  immobile species  via an exchange term in the form 
of a Henry-type law. This is a simple way of modeling a non-equilibrium deposition process. 
 The approach 
 can be extended by using a non-linear differential equation as the one in
\cite{Krehel1}, especially if one aims to come closer to the
experimental data  from \cite{johnson1995dynamics}. Furthermore, each mobile species $u_i$  is transported by diffusion within the patch $Y$ of the porous medium $\Omega $ with a production term expressing the aggregation and fragmentation processes.  
To fix ideas, we assume that  the rate of collision is proportional  to the concentration  of the colliding species and assume the Smoluchowski ansatz.

According to the above description,   the model equations  have the form:
\begin{align}
  \label{eq:model-dim-1} 
  \partial _t u_i&=\Diff_i \Delta u_i + R_i(u)      &&\text{in }Y-B(r)
 % ^\varepsilon
   ,i\inRange{N}\\
  \label{eq:model-dim-2}
  -\Diff_i\nabla u_i\cdot n&=\dca_i u_i-\dcb_i\dep &&\text{on }
  {\Gamma\left(r\right)}\\
  \label{eq:model-dim-3}
  \partial _t \dep&= \sum_{i=1}^N \left(\dca_i u_i-\dcb_i\dep\right) &&\text{on }
  { \Gamma\left(r\right)},
\end{align}
 with $a_i$, $b_i$ constants of proportionality in Henry's law and $d_i$ the diffusion coefficient of the species $i$. 
The reaction rate $R_i(u)$ is here assumed to take the classical   Smoluchowski form (see \cite{Krehel1}) :
\begin{eqnarray}\label{eqreaction}
R_i(u)=\frac12\sum_{i+j=k}^{\infty}\alpha_{i,j}\beta_{i,j}u_i u_j -u_k\sum_{i=1}^\infty \alpha_{k,i}\beta_{k,i}u_i,
\end{eqnarray}
with $u=\left( u_1,u_2,\ldots,u_N\right)$.
 Also $\beta_{i,j}$ is the collision kernel - rate constant determined by the transport mechanism that  bring the particles in close contact, 
$\alpha_{i,j}\in [0,\,1]$ is the collision efficiency, the fraction of collision that finally form an aggregate (\cite{Krehel1}). 
Note that if we consider only one colloidal species, then this reduces to 
$R(u)=-c u^2$, with a suitable reaction constant $c>0$.  As alternative to \eqref{eqreaction}, $R$ can take the Becker-D\"oring form (\cite{Korevaar}).

The growth of the pore radius $r$ relates naturally  to the amount of matter deposited $v$.
 More precisely, the rate of increase of the volume  of immobile species is proportional to the rate of change of deposition of $v$. 
 
Therefore if the volume occupied by the immobile species located at $x\in\Omega$ is denoted with $\mathcal{V}(x,t)$, 
\begin{equation}
{ \partial_t}\mathcal{V}(x,t)=\rho^{-1}\int_{\Gamma(x,t)}{ \partial_t} v,
\end{equation}
for  $\rho$ the density of the immobile species.
 
  Thus for a two-dimensional  radially symmetric setting
  %scenario, 
  $V=\pi r^2$ and   
  for $r'=\frac{d r}{dt}$,
we have
\begin{align}
  \label{eq:model-dim-4}
  r(t,x)r'(t,x)= \rho^{-1}\frac{1}{2\pi}\int_{\Gamma(t,x)}\partial _t\dep.
\end{align}
%Obviously, adopting (\ref{eq:model-dim-4}), we deliberately decide not to include in the discussion eventual surface effects . 
%\rho
%for $r(t,x)$ $\sim$ $r^{\epsilon}(t,x)$}.

The model needs to be completed with initial conditions. At the boundary of the patch $Y$ we consider periodic boundary conditions for $u_i$.  We will discuss in  section \ref{large} the case $\partial Y \cap \partial \Omega\neq \emptyset.$ 

\section{Scaling} In this section, we choose to work with dimensionless quantities and introduce our notation describing the  geometry of the porous material.
 
\subsection{Nondimensionalization at the patch level}\label{nondim}
We set 
$t=t_0\tilde{t}$,
$x=\ell\tilde{x}$,
$\coll_i=\coll_0\tilde{\coll}_i$,
$\dep=\dep_0\tilde{\dep}$.
%\color{blue}
{For convenience, the constants $\Diff_i$, $\dca_i$, $\dcb_i$ can be written in the following way
$\Diff_i=\Diff\tilde{\Diff}_i$,
$\dca_i=\dca \tilde{\dca}_i$,
$\dcb_i=\dcb \tilde{\dcb}_i$, for $d:=\max_{i\in\{1,...,N\}}\{d_i\}, \,a:=\max_{i\in\{1,...,N\}}\{a_i\},\, b:=\max_{i\in\{1,...,N\}}\{b_i\}$}. Also, we take $r=\ell\tilde r$ and, for convenience, we use further $r$ instead of $\tilde r$. The same convention of removing the tilde from $\tilde X$ to yield the scaled quantity $X$ applies also to our scaled geometry (i.e. when referring to $\tilde\Omega$, $\tilde Y$, $\tilde\Gamma$, $\tilde B(r)$,  etc.).

Regarding equation  (\ref{eq:model-dim-1}), dropping tilde in the process, we obtain
\begin{align}
  \label{eq:model-1-sub-1}
  \frac{1}{t_0}\partial _t\coll_i=\frac{\Diff}{\ell^2}\Diff_i\Delta \coll_i+\frac{1}{\coll_0}R_i(\coll_0\coll).
\end{align}
%which is equivalent to
%\begin{align}
%  \label{eq:model-1-sub-2}
%  \partial _t\coll_i=\frac{\Diff t_0}{\ell^2}\Diff_i\Delta \coll_i+\frac{t_0}{\coll_0}R_i(\coll_0\coll).
%\end{align}

By the latter equation we deduce that the characteristic diffusion time scale  is  $t_0=\frac{\ell^2}{d}$.
  However since we want to focus the attention on the deposition process we will make a different choice of  $t_0$ 
in the following. 
 Substituting the indicated scalings into  (\ref{eq:model-dim-2}), we obtain
\begin{align}
  \label{eq:model-2-sub-1}
  -\frac{\Diff}{\ell} \Diff_i\nabla \coll_i\cdot n=\dca(\dca_i\coll_i-\frac{\dcb\dep_0}{\dca\coll_0}\dcb_i\dep).
\end{align}
%and hence,
%\begin{align}
%  \label{eq:model-2-sub-2}
%  - \Diff_i\nabla \coll_i\cdot n=\frac{\dca \ell}{\Diff}(\dca_i\coll_i-\frac{\dcb\dep_0}{\dca\coll_0}\dcb_i\dep).
%\end{align}
Similarly for equation (\ref{eq:model-dim-3}), we get:
\begin{align}
  \label{eq:model-3-sub-1}
  \frac{\dep_0}{t_0}\partial _t\dep
  =
  \coll_0\dca \sum_{i=1}^N
  (\dca_i\coll_i-\frac{\dcb \dep_0}{\dca \coll_0}\dcb_i\dep).
\end{align}
%leading to
%\begin{align}
%  \label{eq:model-3-sub-2}
%  \partial _t\dep
%  =
%  \frac{\dca t_0 \coll_0}{\dep_0}
%  \sum_{i=1}^N
%  (\dca_i\coll_i-\frac{\dcb \dep_0}{\dca \coll_0}\dcb_i\dep).
%\end{align}
Based on  (\ref{eq:model-3-sub-1}), we take as characteristic deposition time scale  $t_0=\frac{\dep_0}{\dca\coll_0}$.
This is the time scale that we choose to non-dimensionalize the equations.

Thus the moving boundary condition (\ref{eq:model-dim-4}) becomes
\begin{align}
  \frac{\ell^2}{t_0} r r' =\frac{\dep_0 \ell}{2\pi t_0 \rho}\int_{\Gamma(r)}\partial _t\dep.
\end{align} %\rho
%Hence, we obtain
%\begin{align}
%  r r' = \frac{\dep_0}{2\pi \ell\rho}\int_{\Gamma(r)}\partial _t\dep.
%\end{align}%\rho
Such a moving boundary condition is sometimes refereed to as {\em kinetic condition}.

Here are the equations grouped together after choosing the deposition time scale:
\begin{align}
  \label{eq:model-inter-1}
  &\partial _t\coll_i=\frac{\Diff \dep_0}{\ell^2\dca\coll_0}\Diff_i\Delta \coll_i+\frac{\dep_0}{\dca \coll_0^2}R_i(\coll_0\coll) \text{ in }Y-B(r),\\
  \label{eq:model-inter-2}
  &- \Diff_i\nabla \coll_i\cdot n=\frac{\dca \ell}{\Diff}(\dca_i\coll_i-\beta_i \dep) \text{ on }
  {\Gamma\left(r\right)}\\
  \label{eq:model-inter-3}
  &\partial _t\dep = \sum_{i=1}^N (\dca_i\coll_i-\beta_i \dep) \text{ on }
  {\Gamma\left(r\right)}\\
  \label{eq:model-inter-4}
  &r r' = \frac{\dep_0}{2\pi \ell\rho}\int_{\Gamma(r)}\partial _t\dep \text{ in }Y-B(r).
\end{align}%\rho

We also denote $\alpha:=\frac{\dep_0}{2\pi \ell}$ %\rho
and $\beta:=\frac{\dcb \dep_0}{\dca \coll_0}\sum_{i=1}^N b_i=\sum_{i=1}^N \beta_i$,
 for $\beta_i:= \frac{\dcb \dep_0}{\dca \coll_0} b_i$. 
To simplify the writing we also change notation $\tilde{R}_i(\coll):=\frac{\dep_0}{\dca \coll_0^2}R_i(\coll_0\coll)$ and finally drop the tilde. 
 Due to the structure of the reaction terms, we have  the dimensional parameters $\frac{v_0}{\dca }\alpha_{i,j}\beta_{i,j}$ and 
$\frac{v_0}{\dca }\alpha_{k,i}\beta_{k,i}$ in front of the  factors $u_i u_j$ and $u_k u_i$, respectively. 
Additionally, we introduce 
$$\kappa:=\frac{\Diff \dep_0}{\ell^2\dca\coll_0}$$ as
the {\em surface Thiele modulus},  i.e. the ratio of the deposition rate (seen as surface reaction) over the diffusion rate; see  \cite{Froment} for more on Thiele-like moduli and related Damk\"ohler numbers.

\subsection{Scaling the geometry} We build the porous material by replicating the information at the patch level $Y$ until covering perfectly $\Omega$\footnote{This is an assumption. In general, pavements with squares can cover perfectly only regions  $\Omega$ that are rectangles. General shapes of $\Omega$ can only be approximately covered, with a certain controllable error if $\partial \Omega$ is sufficiently smooth. }. 
 We assume that  $\varepsilon :=\frac{\dca \ell}{\Diff}$ is a
small positive number. 
%_______________
Indeed we have for example $\epsilon\simeq 7.61e-7$ for values given in  ~\cite{Krehel1}, ( $\ell=0.101\,$m, $\Diff=1.8735e-06\,$m$^2$/sec, $\dca=1.411e-11$m/sec).
%________________________

Note that $\Omega^{\varepsilon}$ is modelled as a composite periodic structure with $\varepsilon>0$ the small scale parameter.
%} We denote by $\varepsilon=\frac{\dca \ell}{\Diff}$. 
We also assume that   the ratio of microscopic and macroscopic length scales, $\frac{d}{l}=O(\varepsilon)$.

  The latter also indicates that  our asymptotics will be connected 
to a parameter regime where $\ell\propto \frac{d}{\sqrt{a}}$ as well as that   $a=\epsilon\frac{d}{l}=O(\varepsilon^2)$.

%Using $\varepsilon=\frac{\dca \ell}{\Diff}$, we additionally take $\frac{d}{\ell}=O(\varepsilon)$ to characterize  the considered structure of the perforated domain 
%(see Figure \ref{fig:locally-periodic}). Note that our asymptotics will be connected 
%to a parameter regime where $\ell\propto \frac{d}{\sqrt{a}}$.

 We consider a periodic array of cells made of  scaled identical patches labeled $Y$. Scaling this patch by $\varepsilon$ defines our standard (homogenization) cell. Translating this standard cell in space covers the porous medium $\Omega$.  For convenience of notation, we write down $r$ as $r^{\epsilon}$ pointing out this way the dependence of $r$ on $\varepsilon$. 

 Let $k=(k_1,k_2)\in\mathbb{Z}^2$ be a vector of indices and  $\{e_1,e_2\}$  be unit vectors in $\mathbb{R}^2$. For $X\subset Y$, we denote by $X^k$ the shifted subset $X^k=X+\sum_{i=1}^2 k_i e_i$.  
   Using the shifted subset notation, the set $\Omega_0^\varepsilon=\cup\{ \varepsilon Y_0^k \, : \, Y_0^k\subset \Omega^\varepsilon \}$ represents the array of pores (voids),
 while $\Gamma^\varepsilon=\partial \Omega_0^\varepsilon $ represents the pores boundary. We describe the skeleton of the porous matrix  in $\Omega$ by  $\Omega^\varepsilon=\Omega\setminus \bar{\Omega_0^\varepsilon}$.
  Also the boundary of the cells of the periodic domain is the set $\Gamma_e=\cup\{ \varepsilon \partial Y^k \}$ .
 In these terms, the colloidal species can be either mobile, $u^{\varepsilon}$, hosted in $\Omega_0^\varepsilon$, or immobile, $v^{\varepsilon}$, sedimented on
$\Gamma^\varepsilon=\Gamma^\varepsilon(r^{\varepsilon})$.

%_____________________

Finally, we obtain the following set of nondimensionalized equations
valid in the locally periodic domain $\Omega^{\varepsilon}$
:
\begin{align}
  \label{eq:model-nodim-1}
  \partial _t u_i^{\varepsilon}&=\kappa d_i \Delta u_i^{\varepsilon} + R_i(u^{\varepsilon})      &&\text{in }\Omega_0^\varepsilon ,i\inRange{N},\\
   \label{eq:model-nodim-2}
  -d_i\nabla u_i^{\varepsilon}\cdot \normal_\varepsilon &=\varepsilon 
  (\dca_i u_i^{\varepsilon}- \beta _i v^{\varepsilon}) &&\text{on }\Gamma^{\varepsilon},\\
  \label{eq:model-nodim-3}
  \partial _t v^{\varepsilon}&= \sum_{i=1}^N \dca_i u_i^{\varepsilon}-\beta  v^{\varepsilon} &&\text{on }\Gamma^{\varepsilon},\\ 
  % \dcb_i
  \label{eq:model-nodim-4}
  & r^{\varepsilon} {r^{\varepsilon}}' = \alpha \int_{\Gamma^{\varepsilon}(r^{\varepsilon})}\partial _t v^{\varepsilon} &&\text{in }\Omega_0^\varepsilon .
\end{align}
The boundary condition at the cell boundary $\Gamma_e$ is
\begin{eqnarray}
-d_i\nabla u_i^{\varepsilon}\cdot \normal_\varepsilon =0,\quad \text{on } \Gamma_e.
\end{eqnarray}
As initial conditions we take

\begin{align}
  \label{initu}
  u_i^{\varepsilon}(0,x)={u_i}^0(x), &&\text{in }\Omega_0^\varepsilon ,i\inRange{N},  \\
   \label{initv}
   v^{\varepsilon}(0,x)={v}^0(x), &&\text{on }\Gamma^{\varepsilon}, \\
  \label{initr}
 r^{\varepsilon}(0)={r}^0.
\end{align}
The precise definition of $\normal_\varepsilon$ will be done in Section \ref{LPH}.  Mind that the sets $\Omega_0^\varepsilon$ and  $\Gamma^{\varepsilon}$ depend explicitly in time via their dependence on $r^\varepsilon(t)$.

Note that the free boundary condition  (\ref{eq:model-nodim-4}) assumes both a small growth speed  and a maximum radius $r$ (i.e. $r<\varepsilon /2$ see cf. \cite{Krehel2}).
We will show in Section \ref{large} how to relax this restriction by using the upscaling technique by  A. Lacey (see \cite{Lacey1}, e.g.).  
Additionally, it is worth mentioning that this $\varepsilon$-dependent free boundary problem (FBP) is expected to be locally weakly solvable (up to the time that the free boundary remains a circle) for any fixed value of the small parameter $\varepsilon>0$;  we refer the reader to Appendix \ref{app} for our arguments.

\section{Locally periodic homogenization}\label{LPH}

\begin{figure}[h!]
  \centering
\includegraphics[bb= 20 20 350 350, scale=.6]{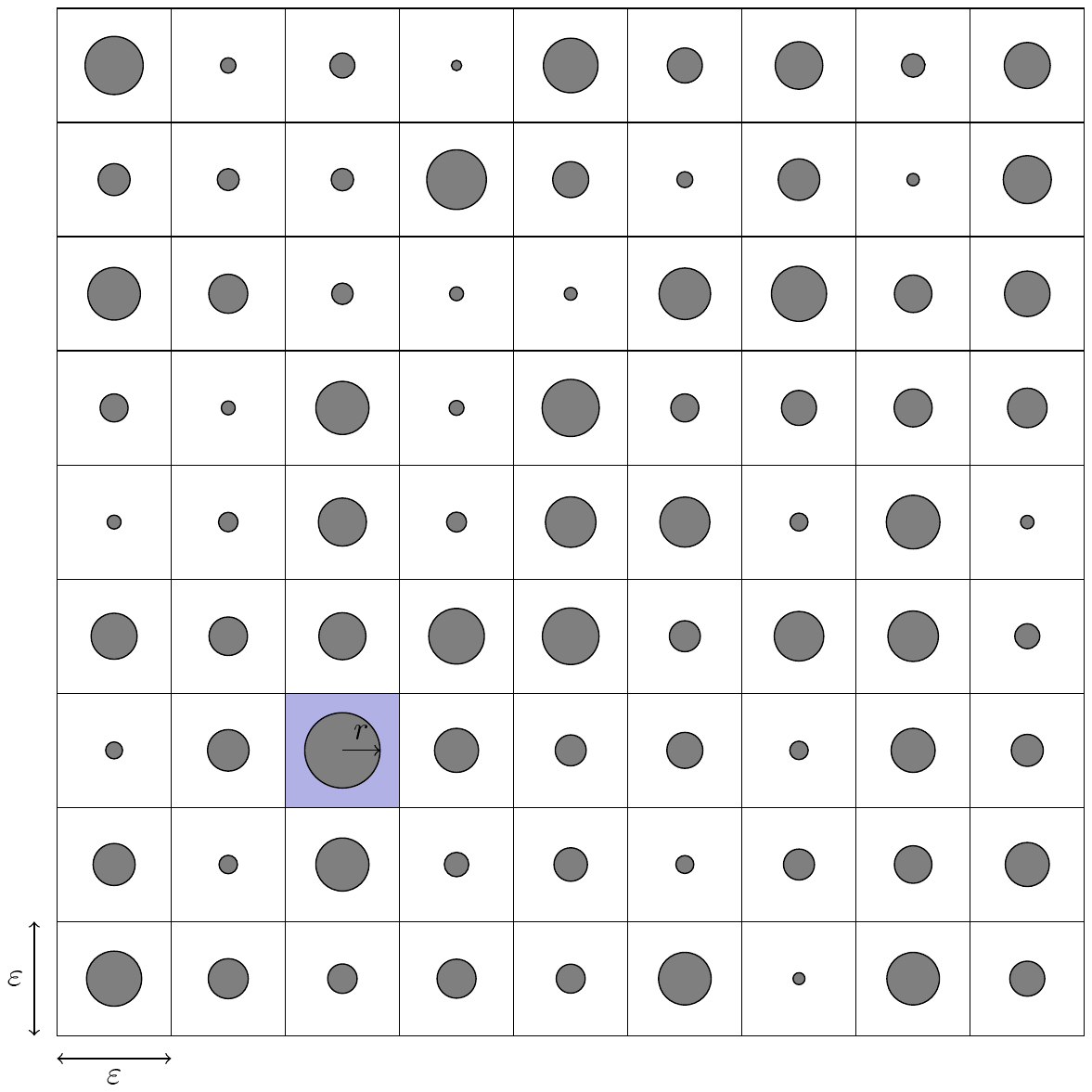}
%bb= 330 230 250 600
 % \includegraphics[width=.4\textwidth]{locally-periodic.ps}     %\vspace*{0.9cm}%locally-periodic.pdf
  \caption{\it Locally periodic pore geometry at time $t=0$. The ``pores" are balls positioned at $x\in \Omega$ with  radius $r(x)>0$. }
  \label{fig:locally-periodic}
\end{figure}

The geometry of our locally periodic structure is described in Figure~\ref{fig:locally-periodic}.
%-----------------------------------------------------------------------------
%\begin{figure}[h!]
%\input{epsf}
%\begin{center}
%\epsfysize=8cm \epsfxsize=8cm \epsfbox{locally-periodic.ps} \vspace*{0.9cm}
%\end{center}
%\caption{\it Locally periodic pore geometry at time $t=0$. The ``pores" are balls positioned at $x\in \Omega$ with  radius $r(x)>0$. {\color{blue}} }\label{fig:locally-periodic}
%\end{figure}
%-----------------------------------------------------------------------------
%____________________________________
%______________________________________
Our upscaling arguments  from this section follow  the lines of \cite{fatima2011homogenization} and \cite{ray2013drug}. 
It is worth pointing out that a rigorous justification of the homogenization asymptotics can be obtained for a given distribution of radii $r(t,x)$ 
(i.e. $r(t,x)$ is fixed and  not part of the problem as in here, in equation \eqref{eq:model-nodim-4})
in terms of the concept of $\theta-2$ convergence of Alexandre (compare \cite{Alexandre}) or using an adaptation of the two-scale convergence designed by G. Nguetseng and G. Allaire to the case of locally periodic functions as in \cite{Ptashnyk} (see references cited therein). Related ideas appear in the context of parabolic equations with rapidly pulsating perforations; see e.g. Ch. 3, Sect. 17 in \cite{CPS}. 

Taking $y=\frac{x}{\varepsilon }$ as the fast variable, we assume the following
formal expansion  (in terms of $Y$-periodic functions)  to hold:
\begin{align}
  \label{eq:formal-expansion}
  u_i^{\varepsilon}(x,y)=\coll_i^0(x,y)+\varepsilon \coll_i^1(x,y)+\varepsilon ^2\coll_i^2(x,y)+\ldots
\end{align}

Noting that $\nabla =\nabla _x+\frac{1}{\varepsilon }\nabla _y$ and substituting the expansion
(\ref{eq:formal-expansion}) into equation (\ref{eq:model-nodim-1}), we obtain:
\newcommand\colls{\coll_i^0+\varepsilon \coll_i^1+\varepsilon ^2\coll_i^2}
\newcommand\Nabla{(\nabla _x+\frac{1}{\varepsilon }\nabla _y)}
\begin{align}
  \label{eq:model-subs-1}
  \partial _t(\colls)=& \kappa \Nabla\cdot
                \left[
                \Diff_i
                \Nabla(\colls)
                \right]
                \nonumber\\
             &+R_i(\coll^0+\varepsilon \coll^1+\varepsilon ^2\coll^2) +\bigO{3}.
\end{align}
% $\nabla d_i \nabla u_i$ or
%____________________________________
{
%\color{red} Shouldn't be from the beginning  $d_i\Delta u_i$\\
%According to \cite{Krehel1} $d_i=\frac{kT}{6\pi n r_i}$, $r_i=i^d r_1$ and $d_i=\frac{1}{i^d} r_1$.
% What $r_i$ represents? If it is the radius of aggregate but 
%mobile species (referring to $u_i$) then it should be a constant and not depending on $y$ later on the analysis in our case (in general it could). If $r_i$ is the radius of the aggregate
%in the center of the cell, i.e. the immobile species then we should have $d_i=d_i(r)$ and a dependence on $\Gamma(x,t)$ but again not in $y$.\\
%I think the first is the right one. Is it so?\\
  Here $u^j=(u_1^j,u_2^j,\cdots,u_N^j)$, $j=0,1,2,\ldots$.
Note also that in our case we have a uniform growing of the perforations, so we are targeting only isotropic situations;
see the geometry of $\Omega_0^\varepsilon$ as well as Figure 1 of  \cite{Krehel1}.
%_____________________________
}

We also apply the same formal expansion for the normal vector,
 pointing out of the surface $\Gamma^{\varepsilon}\left(r(x,t)\right)$, i.e.
\begin{align}
  \label{eq:formal-expansion-normal}
  \normal_\varepsilon(t,x,y)=\normal^0(t,x,y)+\varepsilon \normal^1(t,x,y)+\bigO{2},
\end{align}
where $t$ is just a parameter in the structure of $\normal^0$ and $\normal^1$;  see e.g. \cite{Chechkin} or \cite{fatima2011homogenization}  for the precise analytical expressions of  $\normal^0$ and $\normal^1$ provided 
 that the form 
%the level set
 of the free boundary is known {\em a priori}.

Substituting (\ref{eq:formal-expansion}) and
(\ref{eq:formal-expansion-normal}) into the boundary condition
(\ref{eq:model-nodim-2}) yields:
\begin{align}
  \label{eq:model-subs-2}
  -\Diff_i\Nabla(\colls)\cdot (\normal^0+\varepsilon \normal^1)=&\varepsilon (\dca_i(\colls)-\beta\dcb_i\dep)+\bigO{3}.
\end{align}

By balancing the terms proportional to $\varepsilon ^{-2}$ in (\ref{eq:model-subs-1}) and then those proportional to $\varepsilon ^{-1}$ in (\ref{eq:model-subs-2}), we obtain
\begin{align}
  \label{eq:model-coll0-1}
 {\kappa} \nabla _y\cdot (\Diff_i\nabla \coll_i^0)&=0,\\
  \label{eq:model-coll0-2}
  -\Diff_i\nabla _y\coll_i^0\cdot \normal^0&=0.
\end{align}
Solving (\ref{eq:model-coll0-1})-(\ref{eq:model-coll0-2}) results in 
$u_i^0=u_i^0(x)$,
i.e.  $u_i^0$ is constant in the $y$ variable.

The corresponding terms for $\varepsilon ^{-1}$ in (\ref{eq:model-subs-1}) and $\varepsilon ^{0}$ in (\ref{eq:model-subs-2}) give
\begin{align}
  \label{eq:model-coll1-1}
  {\kappa}\nabla _y\cdot (\Diff_i\nabla _y\coll_i^1)&=0, \\
  -\nabla _y\coll_i^1\cdot \normal^0&=\nabla _x\coll_i^0\cdot \normal^0.
  \label{eq:model-coll1-1b}
\end{align}
 Note that the term $\nabla _y\coll_i^0\cdot \normal^0$ in equation (\ref{eq:model-coll1-1b}) is ommited since $u_i^0=u_i^0(x)$.

 Denote by  $w:=(w_1,w_2)$ the vector of {\it cell functions}. The components  $w_j$ ($j\in\{1,2\}$) are solving the following {\it cell problems}: %\color{red}\color{red}
\begin{align}
  \label{eq:cell-problem-1}
  -\nabla _y\cdot (\Diff_i\nabla _y w_j)&=0,\\
  \label{eq:cell-problem-2}
  \nabla _y w_j\cdot \normal^0&={ e_j}\cdot \normal^0,
\end{align}%\color{red}\color{red}
where $\{e_1,e_2\}$  is the orthonormal basis in $\mathbb{R}^2$.  To pick a concrete unique solution to this Neumann problem, just fix the mean of the solution. 
%{\color{red} It was $e_i$ 
%% as in the paper \cite{Krehel1}.
%I am confused here. $i$ and $d_i$ refers to the mobile colloidal species $u_i$
%$i=1\ldots N$. If it is so  $w_j$, $j=1,2$ correspond to the components of the orthonormal basis  $e_j$, $j=1,2$. Is it so and the index in $e_i$ is a misprint or I have missed something here.}
Then the following representation  holds:
\begin{align}
  \label{eq:model-coll1-2}
  \coll_i^1(x,y)=w(y)\cdot \nabla _x\coll_i^0(x) + \varphi_i(x),
\end{align}
where $\varphi_i$ are given $y$-independent functions. 
Finally, collecting the terms multiplying $\varepsilon ^{0}$ in (\ref{eq:model-subs-1}) and then those multiplying $\varepsilon ^{1}$ in (\ref{eq:model-subs-2}) lead to:
\begin{align}
  \label{eq:model-eps0-1}
  \partial _t\coll_i^0=\Diff_i\Delta _x\coll_i^0+\Diff_i\nabla _x\cdot \nabla _y\coll_i^1+\nabla _y\cdot (\Diff_i(\nabla _x\coll_i^1+\nabla _y\coll_i^2))+R_i(\coll^0),
\end{align}
with 
\begin{align}
 -d_i\nabla _x\coll_i^1\cdot \normal^0&-d_i\nabla _y\coll_i^2\cdot \normal^0 
 -d_i\nabla _x\coll_i^0\cdot \normal^1 -d_i\nabla _y\coll_i^1\cdot \normal^1  
% =\varepsilon 
 =(a_i u_i^0-\beta_i v) 
 &&\text{on }\Gamma^{\varepsilon}\\
 -d_i\nabla _x\coll_i^1\cdot \normal^0&-d_i\nabla _y\coll_i^2\cdot \normal^0
 -d_i\nabla _x\coll_i^0\cdot \normal^1 -d_i\nabla _y\coll_i^1\cdot \normal^1  
 =0  &&\text{on }\Gamma_e,
\end{align}
due to \eqref{eq:model-coll1-1b}.

Integrating (\ref{eq:model-eps0-1}) over $Y(x,t)$, we obtain:
\begin{eqnarray}
  \label{eq:model-eps0-2}
  \int_{Y(x,t)} \partial _t\coll_i^0&=&
  \int_{Y(x,t)}\Diff_i\Delta _x\coll_i^0+
  \int_{Y(x,t)}\Diff_i\nabla _x\cdot \nabla _y\coll_i^1+
  \int_{Y(x,t)} \nabla _y\cdot (\Diff_i(\nabla _x\coll_i^1+\nabla _y\coll_i^2))\nonumber\\
  &+&
  \int_{Y(x,t)}R_i(\coll^0).
\end{eqnarray}
 In the latter equation we have assumed that the domain of quantities such as $\coll_i$, are extended with zero value, inside the solid core of the cell $Y_0$.

The terms in (\ref{eq:model-eps0-2}) that do not depend on $y$ can be
taken out from the integral. Remark the presence of the factor $|Y(x,t)|$. This is in fact intimately linked with the concept of {\em total porosity} of the material (see \cite{Bear}).  Since the area of the periodic cell is 4 (for a cell  $[-1,1]\times [-1,1]$), we have  $|Y(x,t)|=4$ and
$|Y_0(x,t)|=4-\pi r(x,t)^2$. Then the classical concept of porosity is $$\phi(x,t):=\frac{|Y_0(x,t)|}{|\Omega|}=\frac{4-\pi r(x,t)^2}{|\Omega|},$$
where $|\Omega|$ is the measure (area) of $\Omega$.  
 In the following we denote by $A(x,t) := |Y_0(x,t)|$ the active area of the distributed microstructures. 

%This gives
%\begin{eqnarray}
%  A(x,t) \partial _t\coll_i^0&=&A(x,t)\kappa \Diff_i\Delta _x\coll_i^0+
%  \kappa\int_{Y(x,t)} \Diff_i\nabla w:\nabla ^2 \coll_i^0+
%  A(x,t) R_i(\coll^0)\nonumber\\
%  &-&
%  \int_{{\Gamma(x,t)}}(\dca_i\coll_i^0-  \beta_i\dep_0),\label{tic}\\
%  \int_{\Gamma(x,t)}\partial _t \dep_0&=&\int_{\Gamma(x,t)}\sum_{i=1}^N \dca_i\coll_i^0-\beta  \dep_0,\label{bic}\\
%  %\dcb_i\
%  r (x,t)r(x,t)' & = &\alpha \int_{\Gamma(x,t)}\partial _t\dep_0,\\
%  A(x,t)&=&4-\pi r(x,t)^2=\phi(x,t)|\Omega|.\label{eq-rv}
%\end{eqnarray}

 It is worth noting that  dividing (\ref{eq:model-eps0-2})
%  (\ref{tic}) 
  by $|\Omega|$ brings  in as natural factor the (total) volumetric porosity $\phi(x,t)$; the obtained upscaled equation resembles what one expects in standard textbooks on porous media \cite{Bear}.  In Section \ref{large}, since we are interested in capturing macroscopically surface effects due to the motion of the microscopic free boundaries, we prefer to focus on where the area factor $A(x,t)$ enters the upscaled model equations. Note also that regarding the reaction term, $R_i(u^0)$, application of the asymptotic analysis results   to first order terms, lead to an expression of exactly  the same form as \eqref{eqreaction}.
Moreover, since  $Y_0=Y_0(r(x,t))=Y_0(r)$ and
 $ A(x,t) =4-\pi r(x,t)^2=\phi(x,t)|\Omega|$, 
 %by equation \eqref{eq-rv}, 
 $r=r(v)$, we now have the implicit relation  $A=A(v)$.
Looking at  (\ref{eq:model-nodim-3}), after the averaging process (see also equation \eqref{mod1av}) 
%(\ref{bic}),
 we see that we have $\tilde{\dep}:=\int_{\Gamma}\dep_0$ as the accumulated mass of macroscopic immobile species. Interestingly, we observe that in the fast reaction limit $\kappa\to+\infty$ the macroscopic system becomes elliptic.
 More specifically,  we have $\kappa=\frac{\Diff \dep_0}{\ell^2\dca\coll_0}=\frac{1}{\varepsilon}\frac{v_0}{u_0 \ell}$ and in the case that $\frac{v_0}{u_0 \ell}=O(1)$ we obtain the aforementioned limiting case. If this is the case, then  the derivation of the macroscopic equations can be done in  with exactly the same way,  resulting in with the elliptic version of equation \eqref{mod1av},
 %\eqref{bic}, 
 while the rest of the model equations maintain their form.

 This limiting case corresponds in fact to a very fast deposition scenario. 

%___________________________________________________________________________________
%____________________________________________________________________________________

\section{Upscaling extensions for large $r$}\label{large}

Summarizing the results from the previous section,
% and applying also the homogenization process in equations  \eqref{eq:model-nodim-3}, \eqref{eq:model-nodim-4} 
 we obtain the following  equations for the macroscopic scale :
 %__________________________ eq 5.1
%\begin{eqnarray}
%\hspace{-1cm}\frac{\partial u_i^0}{\partial t}(x,t)&=&
%\kappa d_i \Delta_x u_i^0(x, t)+\frac{\kappa}{A(x,t)}\int_{Y(x,t)} d_i\nabla_y w : \nabla_x^2u_i^0+ {\color{blue}R_i(u^0)}\nonumber\\
%&&{\color{blue} -\frac{{L}(x,t)}{A(x,t)} \left(a_i u_i^0(x,t)-\beta_i v_0(x,t)\right),}\label{e1}
%\end{eqnarray}
%_________________________
%{\clr
%\begin{eqnarray}
\bse\label{mod1}
\be
\hspace{-1cm}{ \partial_t} u_i^0(x,t) = D_{ijk}(x,t) \Delta_x u_i^0(x, t)
+ R_i(u^0)\nonumber
\ee
\be
\hspace{3cm} -\frac{{L}(x,t)}{A(x,t)} \left(a_i u_i^0(x,t)-\beta_i v_0(x,t)\right),\label{e1}
\ee
%\end{eqnarray}
describing the diffusion of $u_i$ in the macroscopic domain.  The effective diffusion tensor has the form  
 $$D_{ijk}(x,t)= d_i \phi(x,t)\tau_{jk}(x,t),$$ where the entries
  $$\tau_{jk}(x,t)=\int_{Y(x,t)} \left( \delta_{j,k}+\nabla_{y_j}w_k(z,t) \right) dz,$$ 
   build the tortuosity tensor for all $i=1,\ldots,N$, $j,k=1,2$. In addition 
   \be \label{mod1ha}
{L}(x,t)=\int_{\Gamma(x,t)}ds=2\pi r(x,t),\quad
A(x,t)=\int_{ Y_0(x,t)}dy=4-\pi r^2(x,t),\quad \mbox{(in 2D)}
\ee
%_________________________
\be \label{mod1hb}
R_i(u)=\frac12\sum_{i+j=k}\alpha_{i,j}\beta_{i,j} u_i^0 u_j^0 -
 u_k^0\sum_{i=1}^\infty \alpha_{k,i}\beta_{k,i} u_i^0.
\ee
Here, our cell functions $w:=(w_1(x,y,t),w_2(x,y,t))$, assumed to have constant mean, satisfy
 \be \label{mod1hc}
-\Delta_y w_i=0,\quad i=1,2 \quad \mbox{in}\quad  Y_0(x,t),\\
\ee
\be \label{mod1hd}
 - n_0(x,t)\cdot\nabla_y w_i=0, \quad \mbox{on}\quad \Gamma_e,\quad
 - n_0(x,t)\cdot\nabla_y w_i= n_i(x,t), \quad \mbox{on}\quad \Gamma(x,t).
\ee
with  $\Gamma_e:=\partial Y$ being the boundary of the cell $n_0(x,t)=(n_1(x,t),n_2(x,t))$ is the corresponding normal vector. 
Equation (\ref{e1})
%--(\ref{e4}
need to be complemented with corresponding initial and boundary conditions. In the sequel of this section, we focus the discussion on the case of a one dimensional macroscopic domain, i.e. $x\in[0,1]$. We set
 \be \label{mod1b}
 u_i^0(0,t)=\left\{\begin{array}{cc}
 u_i^b>0 & t\in [0,t_0],\\
 0         & t>t_0,
 \end{array}\right., \quad   {u_i^0}_x (1,t)=0,
 \ee
\be \label{mod1c}
  u_i^0(x,0)=u_i^a(x)\geq 0.
  \ee
  We will discuss in Section  \ref{simulation} a couple of concrete choice of suitable initial and boundary conditions.
   
%_____________________________
Moreover we have
\be\label{mod1av}
{ \partial_t} v_0 (x, t)=  \sum_{i=1}^N \alpha_i u_i^0(x,t)-\beta v_0(x,t),
%\left(\right)
\ee
describing the rate of deposition, 
 with some initial condition
  \be \label{mod1bv}
  v_0(x,0)= v_a(x)\geq 0,
  \ee 
 and  
\be  
 r(x,t)\, { \partial_t} r (x,t)=
 \alpha\left( \sum_{i=1}^N a_i u_i^0(x,t)-\beta  v_0(x,t)\right)
 {L}(x,t), \label{mod1R}
\ee
%where ${L}(x,t):=\int_{\Gamma(x,t)}ds=2\pi r(x,t)$.
together with some initial distribution
\be \label{mod2R}
  r(x,0)=r_a(x)>0, \quad 0<x<1,
\ee
 \ese
%} 
%_________________________________
 
 If perforations do not degenerate, i.e. $0<r(x,t)<1$  so pores do not touch, then it holds
\begin{eqnarray}
{ \partial_t} r (x,t)=
2\pi \alpha\left( \sum_{i=1}^N a_i u_i^0(x,t)-\beta  v_0(x,t)\right)\label{x}\label{e4}.
\end{eqnarray}
Relation (\ref{x}) describes the rate of increase of the bulk of the immobile core inside a cell centred at the macroscopic position $x\in\Omega$ at time $t$.
%_______________________________
% In addition, looking at equation \eqref{e1}, {\color{blue}  we notice that }
 In our case, regarding equation \eqref{e1} for a one-dimensional macroscopic consideration,  $D_{ijk}$ and $\tau_{jk}$  take simpler forms as  for instance,
 $D_{ijk}=D_i=\kappa d_i 
\left( 1+ \frac{1}{A(x,t)}\int_{Y(x,t)} \frac{\partial w_1}{\partial y_1}dy \right)$.

%  in the following equation \eqref{mod1a}: 
% \be \label{mod1a}
%\hspace{-1cm}{\clb \partial_t} u_i^0 (x,t)=
%\kappa d_i 
%\left( 1+ \frac{1}{A(x,t)}\int_{Y(x,t)} \frac{\partial w_1}{\partial y_1}dy \right)\frac{\partial^2 u_i^0}{\partial x^2}(x,t)+ R_i(u^0)
%\nonumber
%\ee
%\be
%\hspace{3cm}-\frac{{L}(x,t)}{A(x,t)} \left(a_i u_i^0(x,t)-\beta_i v_0(x,t)\right),
%\quad 0<x<1,  \quad t\geq 0.
%\ee

The factor $\frac{L(x,t)}{A(x,t)}$  appears in (\ref{e1}) as well as in (\ref{mod1R}) when this is  multiplied
% divided
  by $A(x,t)$. It can be seen as an effective reaction constant incorporating a surface porosity factor. All these considerations are in agreement with  what one would expect while applying directly the macroscopic laws of porous media theory, usually discovered by arguments at REV level involving volume averaging techniques and integral closure relations; see, \cite{Bear}, e.g. 

The system of equations (\ref{mod1}) is assumed to apply for a one-dimensional macroscopic domain (while the related microscopic problem is two-dimensional) with Dirichlet boundary conditions at one side and Neumann boundary conditions at the other side . However, other kind of macroscopic boundary conditions can be also adopted.

In addition, equation \eqref{mod1R} applies for a two-dimensional square cell as far as $0\leq r\leq 1$.
In order to derive an equation for $r$ in the case that we have $r>1$ for a two-dimensional cell or even a three dimensional one, we have to go back to the original law behind the derivation of equation  \eqref{mod1R}. More precisely, we have 
that 
\begin{center}
{\it Rate of increase of the volume  occupied by the immobile species  $\propto$ Rate of change of deposition of the colloidal mass.}
\end{center}
If the volume occupied by the immobile species located at $x\in\Omega$ is denoted with $\mathcal{V}(x,t)$, then  we have 
\begin{equation}
{ \partial_t}\mathcal{V}(x,t)= \alpha\int_{\Gamma(x,t)}{ \partial_t} v,
\end{equation}
 where here the constant $\alpha=\frac{1}{\rho}>0$.

\paragraph{Two-dimensional cell} As $r(x,t)$ reaches $1$, it touches 
$\partial Y(x,t)$.  At this point, we are concerned with the question: \emph{ How  does the moving interface $\Gamma(x,t)$ intersect  the cell boundary $\partial Y(x,t)$?}  One way to proceed is to assume that the sides of the cell 
are always tangent to the  moving boundary  $\Gamma(x,t)$ as this evolves. This is the case of the {\em configuration A }(see \cite{Lacey1}).  Another option is to consider that the moving boundary is always a part of a circle with radius 
$r(x,t)$, $1\leq r(x,t)\leq \sqrt{2}$ and that $\Gamma(x,t)$ intersects with $\partial Y$ at some angle. This is the case of {\em configuration B} (see \cite{Lacey1} for calculation details and general playthrough as developed for the classical Stefan problem).  We refer the reader also to the (remotely resembling) setting described in \cite{Fernandez}.  

\begin{figure}[!htb]\vspace*{4cm}
   \begin{minipage}{0.48\textwidth}
     \centering
  % \special{psfile=FigCorr1a.pdf angle=0 voffset=-285 hoffset=-50 hscale=30 vscale=30} 
   \includegraphics[bb= 400 400 200 400, scale=.3]{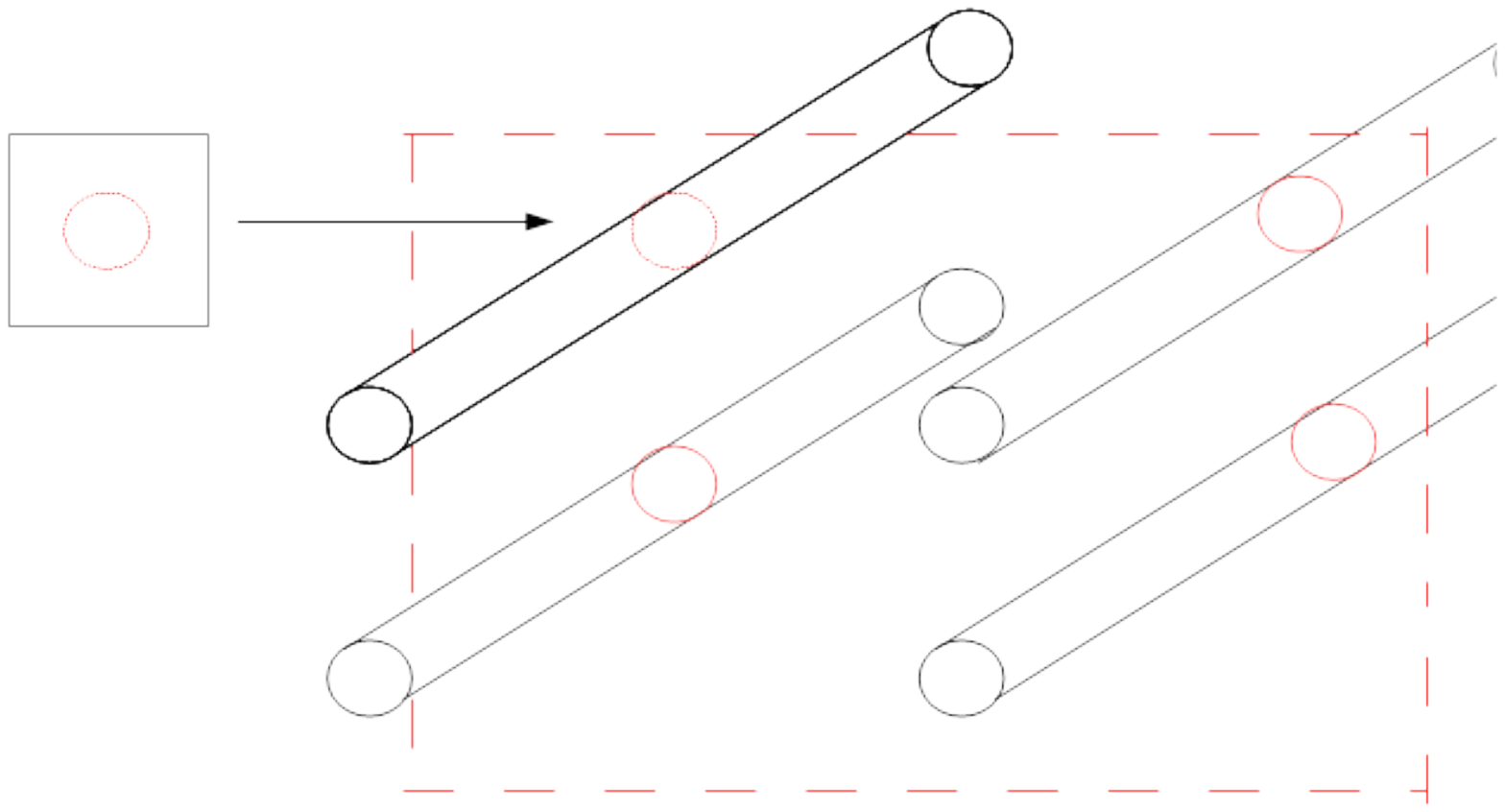}
%\vspace*{-3cm}
\caption{ \it {Schematic representation of a plane crosssecting parallel bars of a solid skeleton. A square on that plane containing the cross section, a circle in this case, centred in it can be considered as one cell in our model. }}\label{cellfig}
   \end{minipage}%\hfill
   \hspace{.5cm}
   \begin{minipage}{0.48\textwidth}
     \centering\hspace{-3cm}
     \includegraphics[bb= 0 350 300 100, scale=.35]{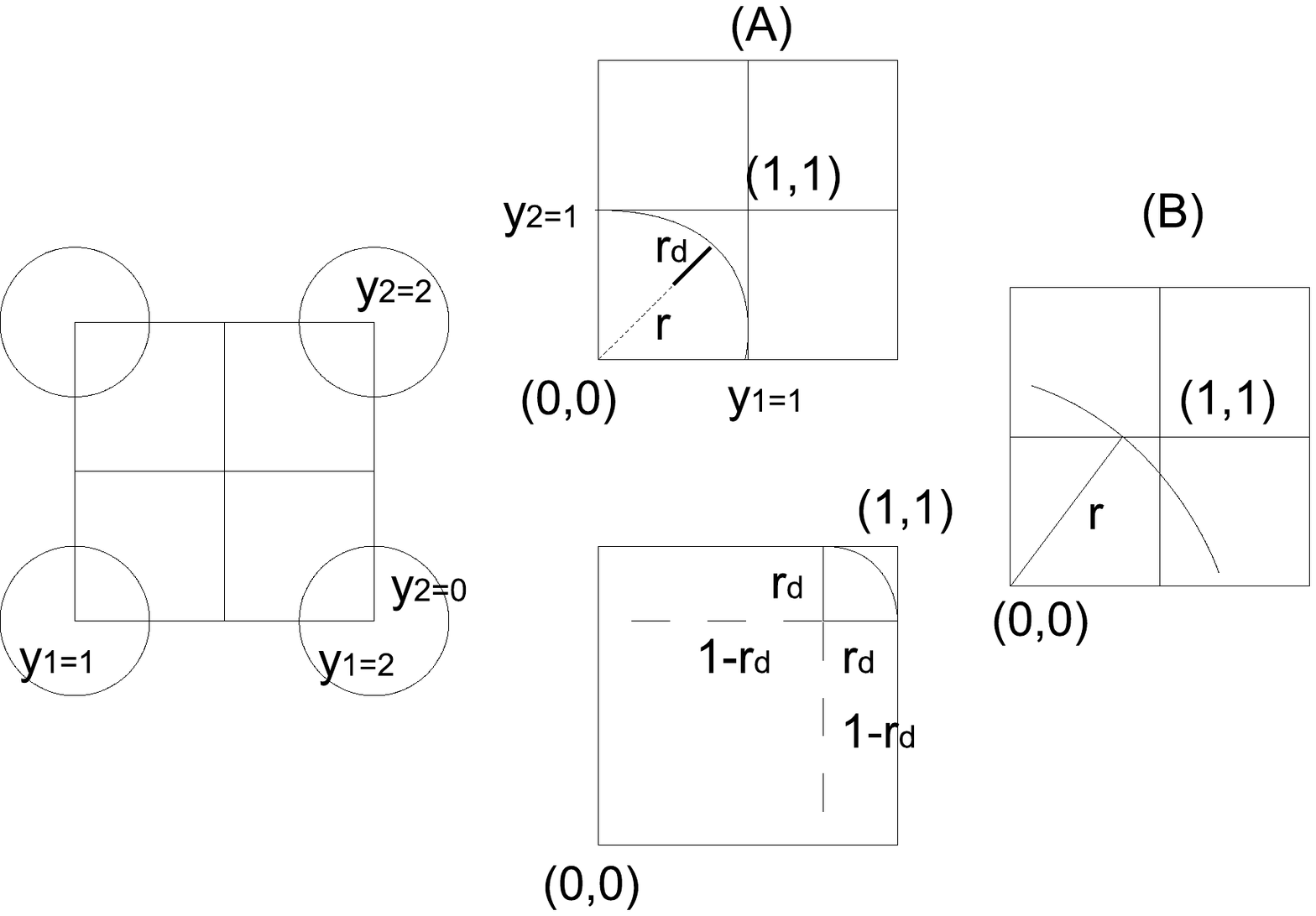}
  \caption{\it Schematic representation of configurations A and B.  In configuration  A we have constant zero angle between the cell boundary and  the  interface while in configuration B the angle of the tangent of the interface and the cell boundary varies as the interface evolves.}%\vspace{-4cm}
 \label{fig6}
   \end{minipage} %\caption{}
\end{figure}

 Moreover, we emphasize that here we assume that clogging takes place when $r$ becomes $\sqrt{2}$. We have in mind that according to this modelling scenario  we still  have diffusion in the direction vertical to the cell plane for $1\leq r\leq \sqrt{2}$. 
This geometrical approach can be justified if we consider inside the material, a solid skeleton consisting of  long narrow cylindrical bars equispaced and parallel (see Figure \ref{cellfig}). 
%__________________________________________________________________
%\begin{figure}[h]%\vspace*{6cm}
%   \includegraphics[bb= 0 350 100 600, scale=.3]{FigCorr1a.pdf}
%%\special{psfile=FigCorr1a.png angle=0 voffset=-185 hoffset=-0 hscale=40 vscale=35} 
%%\vspace*{1cm}
%\caption{ \it {Schematic representation of a plane crosssecting parallel bars of a solid skeleton. A square on that plane containing the cross section, a circle in this case, centred in it can be considered as one cell in our model. }}\label{cellfig}
%\end{figure}
%%corrcell1.png
%__________________________________________________________________

%_____________________________________________________________ 
\vspace{3cm}
%________________________________________________________

%__________________________________
Then we may think of a  plane intersecting these cylinders   transversely  and  obtain a sequence of equispaced circular segments corresponding to the solid parts, while the rest of the plane corresponds to a net of voids.  Next taking advantage of this symmetrical setting we consider   square cells filling the plane with each one of them containing a solid segment centred in it while the rest  is void.
Also, the symmetrical growth of the solid core inside the cell allows us instead of using the standard periodic conditions at the exterior boundary
$\Gamma_e$, to use Neumann conditions as in the first of \eqref{mod1hd}. 
%-----------------------------------------------------------------------------
%\begin{figure}[h]\label{confAB}
%\input{epsf}
%\begin{center}
%\epsfysize=6cm \epsfxsize=12cm \epsfbox{coldyn2.eps} \vspace*{0.9cm}
%\end{center}
%\caption{\it Schematic representation of configurations A and B.}\label{fig6}
%\end{figure}
%-----------------------------------------------------------------------------
%____________________________________
%\begin{figure}[h]
%  \centering
%  \includegraphics[width=.8\textwidth]{coldyn2.eps}\vspace{-7cm}
%  \caption{\it Schematic representation of configurations A and B. {\clr In configuration  A we have constant zero angle between the cell boundary and  the  interface while in configuration B the angle of the tangent of the interface and the cell boundary varies as the interface evolves.}}%\vspace{-4cm}
% \label{fig6}
%\end{figure}
%____________________________________

\paragraph{Configuration A} 
We assume that $1\leq r(x,t)\leq \sqrt{2}$ is the line segment, part of the diagonal  
 of the square (line segment from $(0,0)$ to $(1,1)$ in Figure \ref{fig6}) with end points $(0,0)$ and the point of  intersection of the diagonal with the free boundary $\Gamma(x,t)$. In such a case the free boundary will be a quarter of a circle with radius say $r_d=r_d(x,t)$. Then $r(x,t)$ and $r_d(x,t)$ are related in the following way 
\begin{equation}
(r(x,t)-r_d(x,t))+\sqrt2 r_d(x,t)=\sqrt2.
\end{equation}
Thus 
\begin{eqnarray}
r_d(x,t)=\frac{\sqrt2 -r(x,t)}{\sqrt2 -1},\quad r_d'(x,t)=-\frac{1}{\sqrt2-1}r'(x,t). \nonumber
\end{eqnarray}

In addition, the area occupied by the immobile species will be equal to $4$ times the square quarter minus the square of size $r_d(x,t)$ plus the quarter of the circle with radius $r_d(x,t)$. This gives 
\begin{eqnarray}
\mathcal{V}(x,t)=4\left[1-r_d^2(x,t) +\frac14\pi r_d^2(x,t) \right] = 4\left[1+r_d^2(x,t)(\frac\pi 4 -1) \right].\nonumber
\end{eqnarray}
Therefore the rate of change of the area in terms of $r_d(x,t)$ is 
\begin{eqnarray}
{ \partial_t}\mathcal{V}(x,t)=8\left( \frac \pi 4 -1\right)r_d(x,t) r_d'(x,t)=\alpha \int_{\Gamma(x,t)} v_t.\nonumber
\end{eqnarray}
In terms of $r(x,t)$, the latter relation reads 
\begin{eqnarray}
 \partial_t\mathcal{V}(x,t)=8\frac{\left( 1-\frac \pi 4 \right) }{(\sqrt2-1)^2} (\sqrt2-r(x,t)) r'(x,t) =\alpha \int_{\Gamma(x,t)} v_t,\nonumber
\end{eqnarray}
where $\bar\alpha=\alpha \frac{(\sqrt2-1)^2}{8\left( 1-\frac \pi 4 \right) }$ such that 
\begin{eqnarray}
 (\sqrt2-r(x,t)) r'(x,t) =\bar\alpha \int_{\Gamma(x,t)} v_t\nonumber
\end{eqnarray}
holds. Consequently,  the equation \eqref{mod1R} from the original system needs to be replaced by 
\begin{eqnarray}
 \partial_t r(x,t)=\left\{
\begin{array}{cc}
\frac{1}{  r(x,t) } \alpha\left( \sum_{i=1}^N a_i u_i^0(x,t)-\beta  v_0(x,t)\right){L}(x,t), & r_a\leq r\leq 1,\\
\frac{1}{ \sqrt2- r(x,t) } \bar\alpha\left( \sum_{i=1}^N a_i u_i^0(x,t)-\beta v_0(x,t)\right){L}(x,t), & 1\leq r\leq \sqrt2,
\end{array}\right.
\nonumber
\end{eqnarray}
for 
\begin{eqnarray}
{L}(x,t)=\left\{
\begin{array}{cc}
2\pi r(x,t) , & r_a\leq r\leq 1,\\
2\pi r_d(x,t)=2\pi\frac{\sqrt2-r(x,t)}{\sqrt2-1}, & 1\leq r\leq \sqrt2,
\end{array}\right.
\nonumber
\end{eqnarray}
and 
\begin{eqnarray}
A(x,t)=4-\mathcal{V}(x,t)=\left\{
\begin{array}{cc}
4- \pi r^2(x,t) , & r_a\leq r\leq 1,\\
 4(1-\frac\pi 4 )\left( \frac{\sqrt2 -r(x,t)}{\sqrt2 -1}\right)^2  , & 1\leq r(x,t)\leq \sqrt2.
\end{array}\right.
\nonumber
\end{eqnarray}

\paragraph{Configuration B} 
The free boundary consists here of four arcs of a circle of radius $r(x,t)$, with $1\leq r(x,t)\leq\sqrt2 $. 
In this case the area of the immobile species will be equal to $4$ times the volume of the two triangles plus the circular sector corresponding to angle $\frac\pi 2 -2\cos^{-1}(\frac1r)$. More specifically,  we have 
\begin{eqnarray}
\mathcal{V}(x,t)=4\left[\sqrt{r(x,t)^2-1} +\frac{r(x,t)^2}{2}\left(\frac\pi 2 -2\cos^{-1}(\frac1{r(x,t)}) \right)\right]. \nonumber
\end{eqnarray}

Then the rate of change of the immobile species area will be 
\begin{eqnarray}
 \partial_t\mathcal{V}(x,t) r(x,t)= 
 \frac{d}{dt}\left[4\sqrt{r(x,t)^2-1} +r(x,t)^2\left(\pi  -4\cos^{-1}(\frac1{r(x,t)}) \right) \right]=
 \alpha \int_{\Gamma(x,t)} v_t (x,t),\nonumber
\end{eqnarray}
or
\begin{eqnarray}
2r(x,t)\left(\pi  -4\cos^{-1}(\frac1{r(x,t)}) \right)r'(x,t) =
 \alpha \int_{\Gamma(x,t)} v_t.\nonumber
\end{eqnarray}

In such a case 

\begin{eqnarray}
 \partial_t r(x,t)=\left\{
\begin{array}{cc}
\frac{1}{  r(x,t) } \alpha\left( \sum_{i=1}^N a_i u_i^0(x,t)-\beta  v_0(x,t)\right){L}(x,t), & r_a\leq r(x,t)\leq 1,\\
 \gamma(x,t)\left( \sum_{i=1}^N a_i u_i^0(x,t)-\beta v_0(x,t)\right)L(x,t), & 1\leq r(x,t)\leq \sqrt2,
\end{array}\right.
\nonumber
\end{eqnarray}
with
\begin{eqnarray}
\gamma(x,t)&:=&\frac{\alpha}{ 2 r(x,t)\left(\pi  -4\cos^{-1}\left(\frac1{r(x,t)}\right) \right)}\nonumber\\
{L}(x,t)&=&\left\{
\begin{array}{cc}
2\pi r(x,t) , & r_a\leq r(x,t)\leq 1,\\
r(x,t)\left[2\pi-8\cos^{-1}\left(\frac1{r(x,t)}\right) \right], & 1\leq r(x,t)\leq \sqrt2,
\end{array}\right.
\nonumber
\end{eqnarray}
and, finally, 
\begin{eqnarray}
A(x,t)&=&4-\mathcal{V}(x,t)\nonumber\\
&=&\left\{
\begin{array}{cc}
4-\pi r^2(x,t) , & r_a\leq r\leq 1,\\
4\left[1-\sqrt{r^2(x,t)-1} +r^2(x,t)\left(\frac{\pi}4  -\cos^{-1}(\frac1{r(x,t)}) \right)\right], & 1\leq r(x,t)\leq \sqrt2.
\end{array}\right.
\nonumber
\end{eqnarray}

%_________________________________________________________________________________
%__________________________________________________________________________________

\section{Numerical approximation}\label{simulation}

To treat numerically  problem \eqref{mod1}, we need firstly 
to  obtain a good  numerical approximation for the cell problems  \eqref{mod1hc} and determine the shape of the cell functions  $w_1,w_2$ posed in $ Y_0(x,t)$. 
More specifically, we proceed for the various values of 
$r$, at a first stage for $r_a\leq r(x,t) \leq 1$ and, at a second stage, for $1<r(x,t)\leq \sqrt{2}$ depending which  microscopic configuration (A or B) is considered. We  take a partition of width $\delta r$, $r_a=r_0, r_1=r_0+\delta r,\ldots, r_{M_1}=1$. Then since $ Y_0$ is determined as the area contained inside the square cell and outside  the circle of radius $r$, we obtain a sequence of solutions for each 
$ {Y_0}_i$ corresponding to the radius $r_i$ of the partition. In the same way, for a partition of the interval $[1, \sqrt{2}]$,  $1=r_{M_1+1}, \ldots, r_M=\sqrt{2}$ we solve again the problem in the domain the way that this is specified in the case of configuration A or B. Note that in these cases the domain $ Y_0$  is reduced into four separate segments and thus we need the solution (due to symmetry) in only one of them.
Then  we  use a finite element scheme to solve these cell problems. The finite element numerics  package in MATLAB ''Distmesh" (\cite{Persson}) is used to triangulate the domain ${ Y_0}_i= Y_0(r_i)$ and a solver that has been implemented for this specific problem (equations \eqref{mod1hc}) is applied. 

\begin{figure}[h]
\begin{center}
\includegraphics[bb= 0 0 450 300, scale=.7]{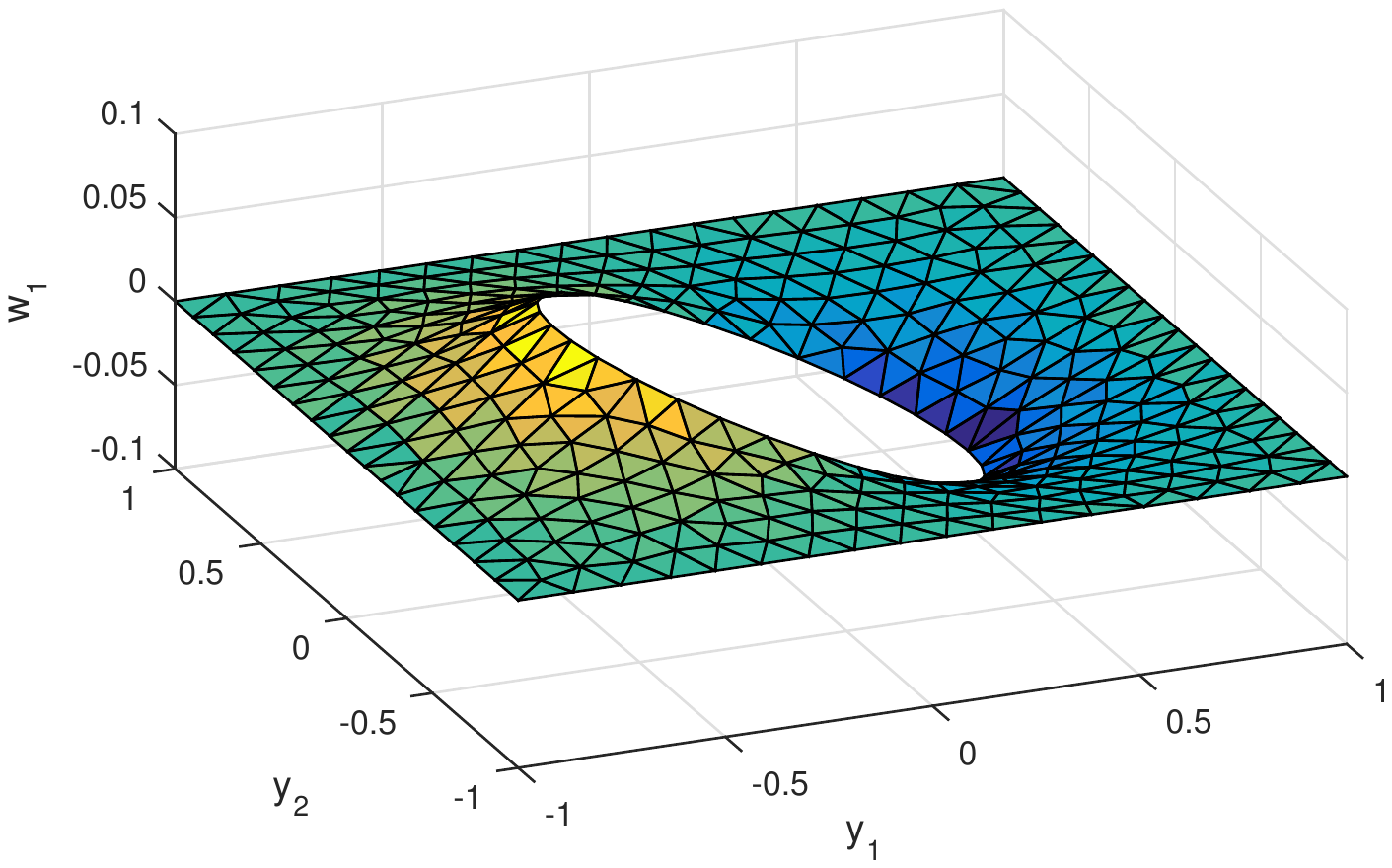}
\end{center}
\caption{\it Numerical solution of the cell problem \eqref{mod1hc} and specifically for $w_1$ with $r_i=.5$.}\label{Figw1}
\end{figure}

\begin{figure}[h]
\input{epsf}
\begin{center}
\includegraphics[bb= 0 0 450 250, scale=.7]{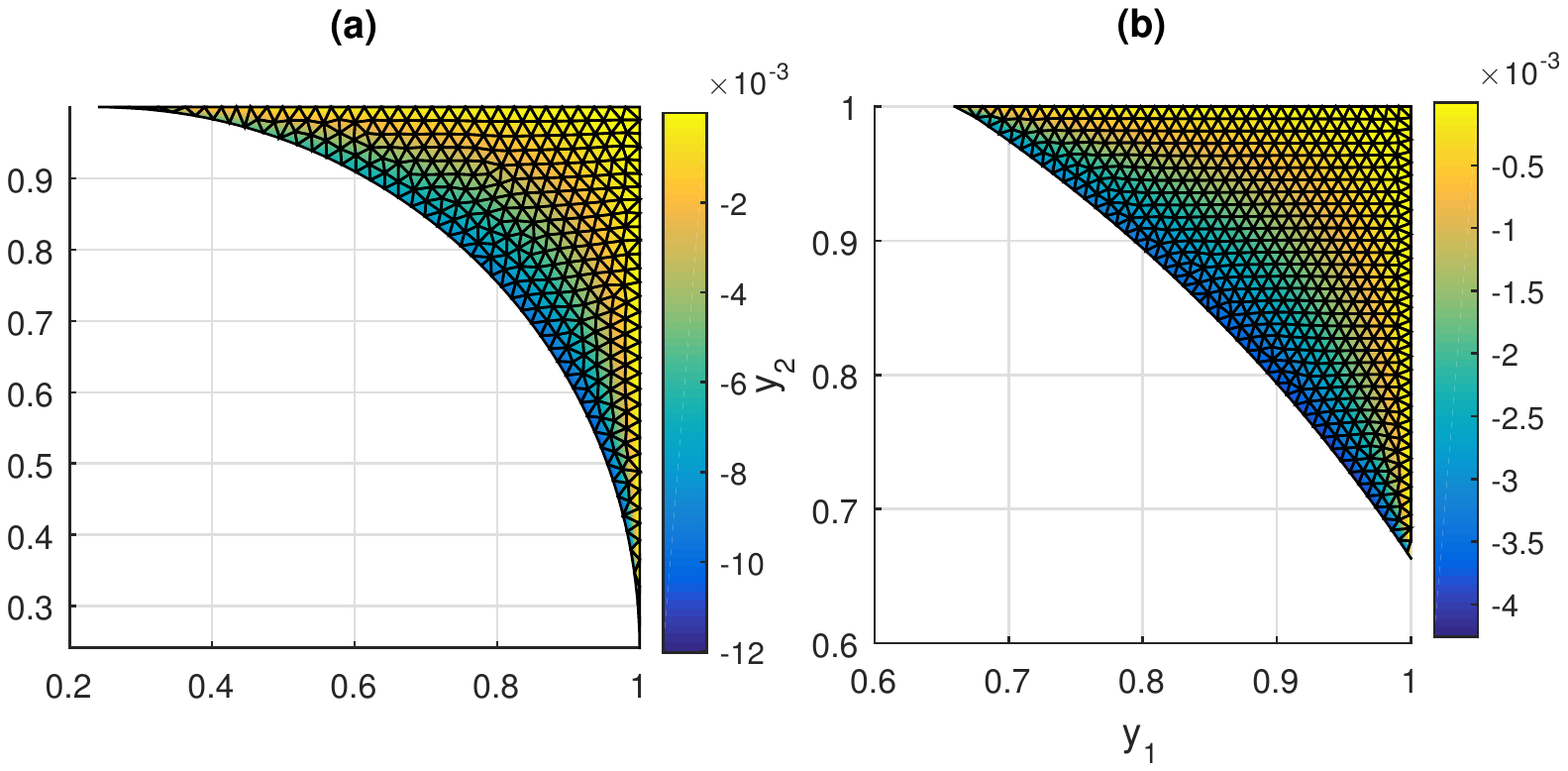}
\vspace*{-0.9cm}
\end{center}
\caption{\it Numerical solution of the cell problem \eqref{mod1hc}: The profile of  $w_1$ for (a) the configuration A and for (b) the configuration B.}\label{Figw2}
\end{figure}

We illustrate the numerical solution for this problem for specific choices of $r_i$.; see  Figure \ref{Figw1} for  $r_i=.5$ and  Figure \ref{Figw2} for $r_i=1.2$ and configuration A, and for $r_i=1.1$ and configuration B.
%-----------------------------------------------------------------------------
%----------------------------------------------------------------------
%____________________________________HHH
%\begin{figure}[h]
%  \centering
%  \includegraphics[width=.9\textwidth]{CPcircle.eps}
%  \caption{\it Numerical solution of the cell problem \eqref{mod1hc} and specifically for $w_1$ with $r_i=.5$.}
%\label{Figw1}
%\end{figure}
%____________________________________
% %_____________________________________________________________ 
%\begin{figure}[!htb]
%   \begin{minipage}{0.48\textwidth}
%     \centering
%   %  \includegraphics[width=6cm,height=6cm]{CPcircle.eps}%\vspace{-6cm}
%     \caption{\it Numerical solution of the cell problem \eqref{mod1hc} and specifically for $w_1$ with $r_i=.5$.}
%     %\label{Figw1}
%   \end{minipage}%\hfill
%   \hspace{.5cm}
%   \begin{minipage}{0.48\textwidth}
%     \centering%\vspace{-3cm}
%%     \includegraphics[width=6cm,height=6cm]{diff2n.eps}%\vspace{-3cm}  
%          \caption{\it The diffusion matrix dependence on the choice of the pore radius.}
%          \label{fig:diffusion-coeff}
%   \end{minipage} %\caption{}
%\end{figure}
%%________________________________________________________
%-----------------------------------------------------------------------------
%----------------------------------------------------------------------
%____________________________________
%\begin{figure}[h]
%  \centering
% % \includegraphics[width=.7\textwidth]{CPconfABex.eps}
%  \caption{\it Numerical solution of the cell problem \eqref{mod1hc}: The profile of  $w_1$ for (a) the configuration A and for (b) the configuration B.}
%\label{Figw2}
%\end{figure}
%____________________________________

\begin{figure}[h!]
  \centering
  \includegraphics[bb= 0 0 350 220, scale=.7]{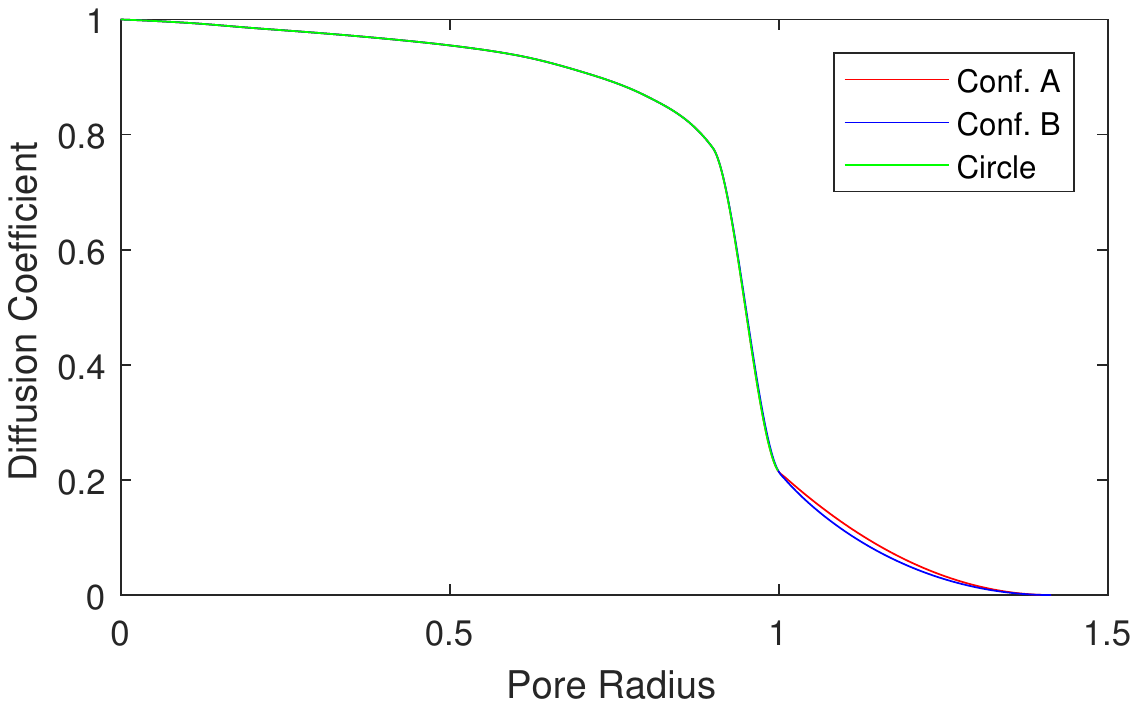}
  \caption{The diffusion matrix dependence on the choice of the pore radius.}
  \label{fig:diffusion-coeff}
\end{figure}

Having available the numerical evaluation of the cell functions $w$ as approximate solutions to  the cell problems \eqref{mod1hc} and \eqref{mod1hd}, the entries of the diffusion tensor  $D_{ijk}=\int_{ Y_0(x,t)} d_i \left( \delta_{j,k}+\nabla_{y_j}w_k \right)$, $i=1,\ldots,N$, $j,k=1,2$ can be calculated directly. More specifically, for each $(x,t)$ and consequently for each $ Y_0(x,t)$, the corresponding value of $D_{ijk}(x,t)$ is approximated via linear interpolation.
For the cases that $0<r(x,t)<1$, ($4-\pi\leq | Y_0(x,t)|\leq 4$)
 as well as for $1<r(x,t)<\sqrt2$ and configurations A and B, the results regarding the diagonal entries, needed for equation 
 %\eqref{mod1a}
  \eqref{e1}, are shown in Figure \ref{fig:diffusion-coeff}.
%__________________________________________________________HHH
%__________________________________________________________

Next, we solve the system of equations \eqref{e1}-\eqref{mod2R}. We use a finite element scheme to solve the  one-dimensional version of the field equation \eqref{e1},
% i.e. equation \eqref{mod1a}
 together with its boundary and initial conditions. 
Let $\psi_j=\psi_j(r)$, $j=0,\ldots, M$ denote the standard linear B - splines defined on the interval $[0,\,1]$ with respect to the considered partition.

 We then take the Galerkin approximation $u_i(x,t)=\sum_{j=0}^M {a}_{{u_i}_j}(t) \psi_j(x)$, 
 $v(x,t)=\sum_{j=0}^M {a}_{{v}_j}(t) \psi_j(x)$
$\tau\geq 0$, $0\leq r(x,t)\leq 1$. We also take  $r(x,t)=\sum_{j=0}^M {a}_{{r}_j}(t) \psi_j(x)$ for $r$ in \eqref{mod1R}.
Applying the standard Galerkin method, we obtain a system of equations for the coefficients ${a}_{{u_i}_j}(t), a_{{v}_j}(t), {a}_{{r}_j}(t)$, $i=1,\ldots,N$, $j=1,\ldots,M$. The resulting system of ODE's 
is finally solved by an implicit time-stepping  scheme.
%____________________________________________________________________
 \paragraph*{\bf Finite element scheme for the model equations}

We substitute the above expressions  of $u_i(x,t)$, $v(x,t)$ and 
$r(x,t)$  in equation 
 \eqref{mod1} 
%\eqref{mod1a} 
multiply with $\psi_l$ and then integrate over $[0,\,1]$ to obtain 

  \begin{eqnarray}\hspace{-.1cm}\label{eqfe_u}
  \sum_{j=0}^M {\dot{a}_{{u_i}_j}}\int_0^1\psi_j\,\psi_{l}\,dx
 &=& \kappa d_i 
 \sum_{j=0}^M a_{{u_i}_j}\int_0^1 D(x,t)\psi'_j\,\psi'_{l}\,dx  
 + \int_0^1 R  \left(\sum_{j=0}^M {a}_{{u_i}_j} \psi_j \right)\psi_{l}\,dx
 \nonumber\\
   && {-}\left(\alpha_i a_{{u_i}_j}-\beta_i a_{{v}_j}\right) 
   \int_0^1 \frac{L(x,t)}{A(x,t)}\psi_j\,\psi_{l} dx
  \end{eqnarray} 
   where  ${l}=1,2,\ldots,M$  and the dot, $``\, \cdot{} \,"$, denotes differentiation with respect to time.
     
   Similarly, for equation \eqref{mod1av}, we have     
     \begin{eqnarray}\hspace{-.1cm}\label{eqfe_v}
     \sum_{j=0}^M {\dot{a}_{{v}_j}}\int_0^1\psi_j\,\psi_{l}\,dx =
     \sum_{i=1}^N \alpha_i 
     \sum_{j=0}^M a_{{u_i}_j}\int_0^1 \psi_j\,\psi_{l}\,dx  
     -\beta \sum_{j=0}^M a_{{v}_j}\int_0^1 \psi_j\,\psi_{l}\,dx. \quad
      \end{eqnarray}       
      % \left(\right)
 Finally, for equation \eqref{mod1R} it holds
  \begin{eqnarray}\hspace{-.1cm}\label{eqfe_R}
  \sum_{j=0}^M {\dot{a}_{{r}_j}}\int_0^1\psi_j\,\psi_{l}\,dx =
      2\pi\alpha\left( \sum_{i=1}^N \alpha_i 
     \sum_{j=0}^M a_{{u_i}_j}\int_0^1 \psi_j\,\psi_{l}\,dx  
     -\beta\sum_{j=0}^M a_{{v}_j}\int_0^1 \psi_j\,\psi_{l}\,dx \right).
     \quad\quad\,
      \end{eqnarray} 

  Setting $a_{u_i}:=[a_{{u_i}_1},\,a_{{u_i}_2},\ldots,a_{{u_i}_M}]^T$ and 
   $a_{v}:=[a_{{v}_1},\,a_{{v}_2},\ldots,a_{{v}_M}]^T$, 
the system of equations for the coefficients ${a_{u_i}}$'s takes then the form
\begin{eqnarray}
B_l {\dot{a}_{u_i}}(t)=  B_r(t)a_{u_i}(t) +b_R(t) {-} 
C(t)\left(\alpha a_{u_i}(t) - \beta a_{v}(t)\right), \nonumber
\end{eqnarray}
with $t\in (0,T]$.

The corresponding  FEM matrices are 

\[
B_l:=\left(\int_0^1\psi_j(x)\,\psi_{l}(x)\,dx \right), \,\,
B_r:=\kappa d_i\left(\int_0^1 D(x,t)\psi'_j(x)\,\psi'_{l}(x)\,dx\right),\]
\[ C:=\left(\int_0^1 \frac{L(x,t)}{A(x,t)}\psi_j(x)\,\psi_{l}(x)\,dx\right),\]
 while $b_R(t) $  is the array with entries 
$b_R(t):= \int_0^1 R  \left(\sum_{j=0}^M {a}_{{u_i}_j} \psi_j \right)\psi_{l}\,dx$.

What concerns the equation of $v$  we have 
\begin{eqnarray}
B_l {\dot{a}_{v}}(t)=  \left( \sum_{i=1}^N \alpha_i B_l
\left(\alpha a_{u_i}(t) - \beta a_{v}(t)\right) \right), \nonumber
\end{eqnarray}
while for $r$ we get 
\begin{eqnarray}
B_l {\dot{a}_{r}}(t)= 2\pi \left( \sum_{i=1}^N \alpha_i B_l
\left(\alpha a_{u_i}(t) - \beta a_{v}(t)\right)\right). \nonumber
\end{eqnarray}

Next, we apply an implicit in time scheme of the form
 \begin{eqnarray}\label{numscha1}
&\left[B_l-\d t\,B_r(t_n)+\d t\,\alpha C(t_n) \right]a_{u_i}^{n+1}
-\d t\,\beta C(t_n)a_{v}^{n+1} =B_l a_{u_i}^n +\d t\, b_R(t_n),\\
& a_{u_i}^{n+1} -\d t\, \left( \sum_{i=1}^N \alpha_i \left(\alpha a_{u_i}^{n+1}
 - \beta a_{v}^{n+1}\right) \right)= a_{u_i}^n ,\\
& a_{r}^{n+1} -\d t\, 2\pi\left( \sum_{i=1}^N \alpha_i 
\left(\alpha a_{u_i}^{n+1} - \beta a_{v}^{n+1}\right) \right)= a_{r}^n, 
\end{eqnarray}
for   $\d t$ the time step and $t_n=(n-1)\d t\in [0, T]$, $n\in\{1,\ldots, [T/n] \}$, to solve the system numerically.

\section{Discussion of simulation results}\label{discussion}

%\paragraph{\bf Numerical illustrations} 

To obtain a couple of relevant simulation examples, we choose the following set of parameter values:
We consider $N=3$  mobile species $u_i$ and one immobile species. 

We take
$\kappa=1,\,\, 
(d_1,\,d_2,\,d_3)=(.3,.5,.99)$,\,\,
$(a_1,\,a_2,\,a_3)=(.9,.5,.3)$,\,\,
$(\beta_1,\,\beta_2,\,\beta_3)=(1,1,1)$,
$(u_1^b,\,u_2^b,\,u_3^b)=(1,1,1)$,\,\,  
$\alpha_{i,j}=.1,\,\,\beta_{i,j}=100$, $i,j=1,\ldots 3$, $u_a^i(x)=0,\,\,v_a(x)=0,\,\, r_a(x)=.1 ,\,\,
0\leq x\leq 1$.  
 Finally, the generation of mass at the boundary $x=0$ is traced  up to  time $t_0=2$.  At this particular instant, we also observe a relevant jump in the form of the $u_i$'s (cf.  the next simulations).
 For all these simulations,  the reference parameters are chosen in the range indicated in ~\cite{Krehel1}.% and also attempting  to emphasize in the variations of the model behaviour.

\subsection{Basic output. Detecting clogging regions}

\begin{figure}[h!]
  \centering
  \includegraphics[bb= 0 0 480 320, scale=.6]{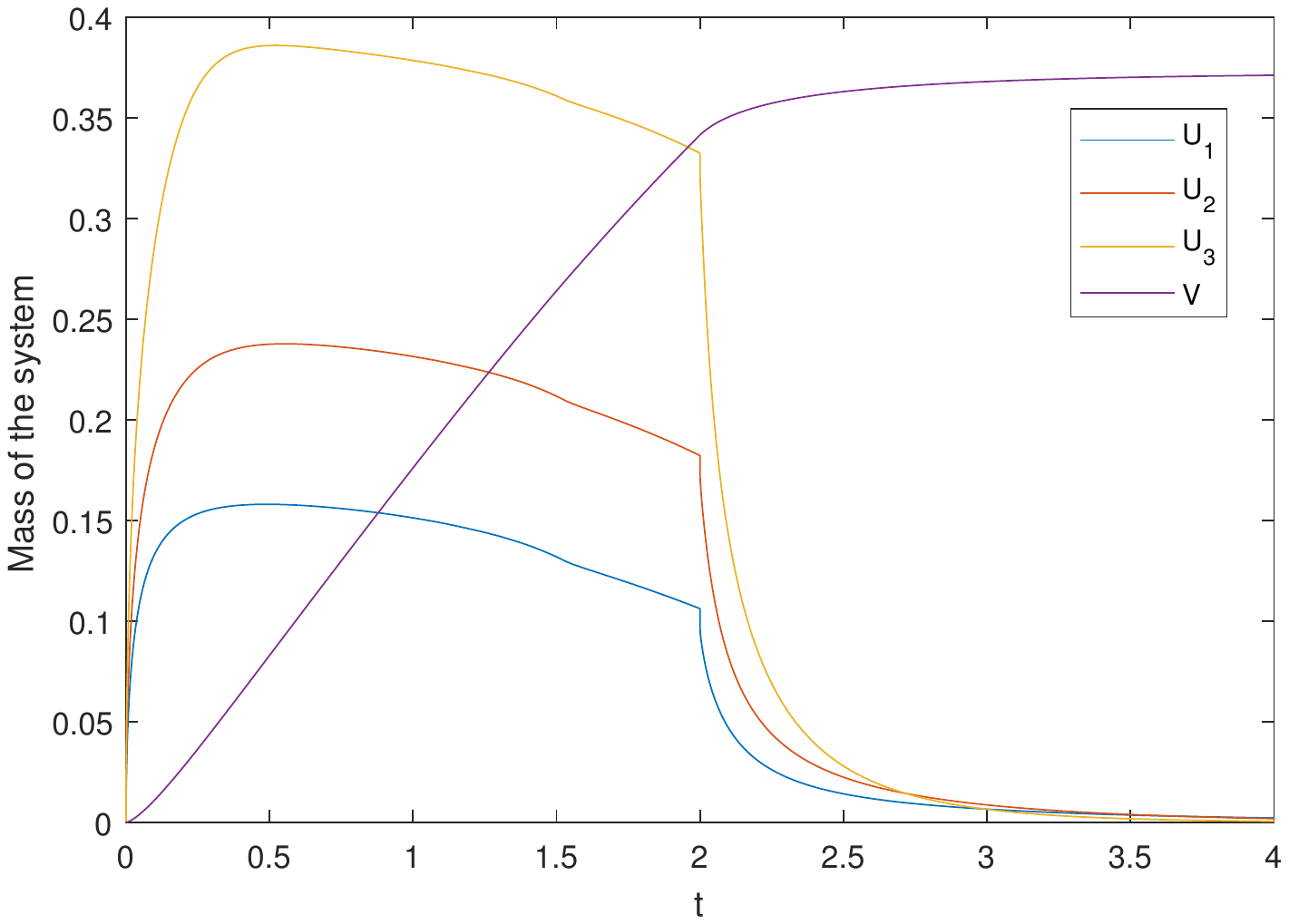}
  \caption{\it Mass of the system against time for the various colloidal species in the system}
\label{Figsim1}
\end{figure}

The basic simulation output is shown in Figure \ref{Figsim1}. This  is done for configuration A, where the mass of the various components of the system, i.e. $U_i(t)=\int_0^1u_i(x,t)dx$, $V(t)=\int_0^1 v(x,t)dx$,  is plotted against time. 

%-----------------------------------------------------------------------------
%\begin{figure}[h!]
%\input{epsf}
%\begin{center}
%\epsfysize=8cm \epsfxsize=12cm \epsfbox{Figsim1a.eps} %\vspace*{0.9cm}
%\end{center}
%\caption{\it Mass of the system against time for the various colloidal species in the system.}\label{Figsim1}
%\end{figure}
%----------------------------------------------------------------------
%____________________________________HHH
\begin{figure}[h!]
  \centering
  \includegraphics[bb= 0 0 480 320, scale=.6]{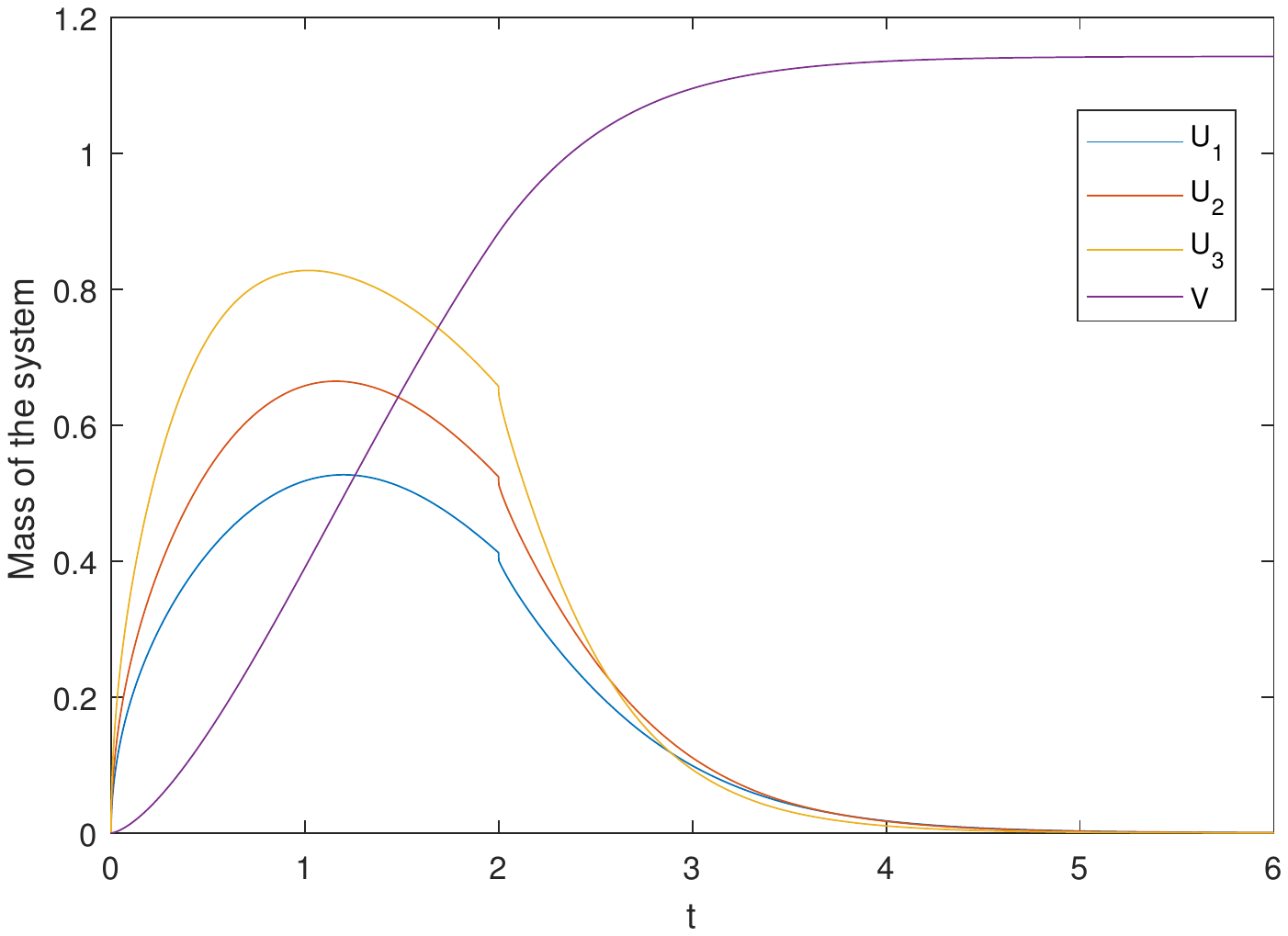}
  \caption{\it Mass of the system against time for the various colloidal species in the system.}
\label{Figsim2}
\end{figure}

%____________________________________
 %_____________________________________________________________ 
%\begin{figure}[!htb]
%   \begin{minipage}{0.48\textwidth}
%     \centering
%  %   \includegraphics[width=6cm,height=6cm]{Figsim1a.eps}%\vspace{-6cm}
%     \caption{\it  Mass of the system against time for the various colloidal species in the system}
%     \label{Figsim1}
%   \end{minipage}%\hfill
%   \hspace{.5cm}
%   \begin{minipage}{0.48\textwidth}
%     \centering%\vspace{-3cm}
% %    \includegraphics[width=6cm,height=6cm]{Figsim2a.eps}%\vspace{-3cm}  
%          \caption{\it Mass of the system against time for the various colloidal species in the system.}
%          \label{Figsim2}
%   \end{minipage} %\caption{}
%\end{figure}
%________________________________________________________
 If we decrease the value of the parameters $\beta_{i,j}$,
  the collision kernel expressing the ability of the particles to aggregate,
and for instance take them to be uniformly constant $\beta_{i,j}=1$, then (after the time that $u_i(0,t)$ becomes zero) the decay of the $U_i$'s is 
%much less steeper  
 moderate
 due to the fact that the aggregation mechanism now is weaker. 
 This is shown in Figure \ref{Figsim2}.
%-----------------------------------------------------------------------------
%\begin{figure}[h!]
%\input{epsf}
%\begin{center}
%\epsfysize=8cm \epsfxsize=12cm \epsfbox{Figsim2a.eps} %\vspace*{0.9cm}
%\end{center}
%\caption{\it Mass of the system against time for the various colloidal species in the system.}\label{Figsim2}
%\end{figure}
%------------------------------------------------------------------------------
%____________________________________HHH
%____________________________________

Increasing the number of mobile species yields similar results. As an example,  for $N=5$
mobile species the simulation is demonstrated in the following Figure \ref{Figsim5}. The parameter values  used here are
the same as those used in the previous simulation in Figure \ref{Figsim1} for the first three species, $i=1,2,3$.
%different than those employed to obtain Figure  \ref{Figsim1}.
 More specifically,  we take 
%\begin{center}
$(d_1,\,d_2,\,d_3,d_4,\,d_5)=(.3,.5,.99,.2,.3)$, 
$(a_1,\,a_2,\,a_3,\,a_4,\,a_5)=(.9,.5,.3,.4,.3)$, 
$(\beta_1,\,\beta_2,\,\beta_3,\,\beta_4,\,\beta_5)=(1,1,1,1,1)$,
 $(u_1^b,\, u_2^b,\, u_3^b,\,$ $ u_4^b,\, u_5^b)=(1,1,1,1,1)$.
  The proximity of the variables $u_1$ and $u_4$, and also of $u_2$ and $u_5$  in Figure \ref{Figsim5}  is due to the particular (arbitrary) choice of the diffusion coefficients $d_i$'s and the $a$'s.
%\end{center}

%$v_a(x)=0$, $r_a(x)=.1$, $0\leq x\leq 1$,\\
%$\alpha_{i,j}=.1$, $\beta_{i,j}=100$.
%-----------------------------------------------------------------------------
%\begin{figure}[h!]
%\input{epsf}
%\begin{center}
%\epsfysize=8cm \epsfxsize=12cm \epsfbox{simfig5sA1.eps} %\vspace*{0.9cm}
%\end{center}
%\caption{\it Mass of the system against time for the various components of the system in the case of 5 mobile species}\label{Figsim5}
%\end{figure}
%----------------------------------------------------------------------
%____________________________________
\begin{figure}[h!]
  \centering
    \includegraphics[bb= 0 0 600 420, scale=.47]{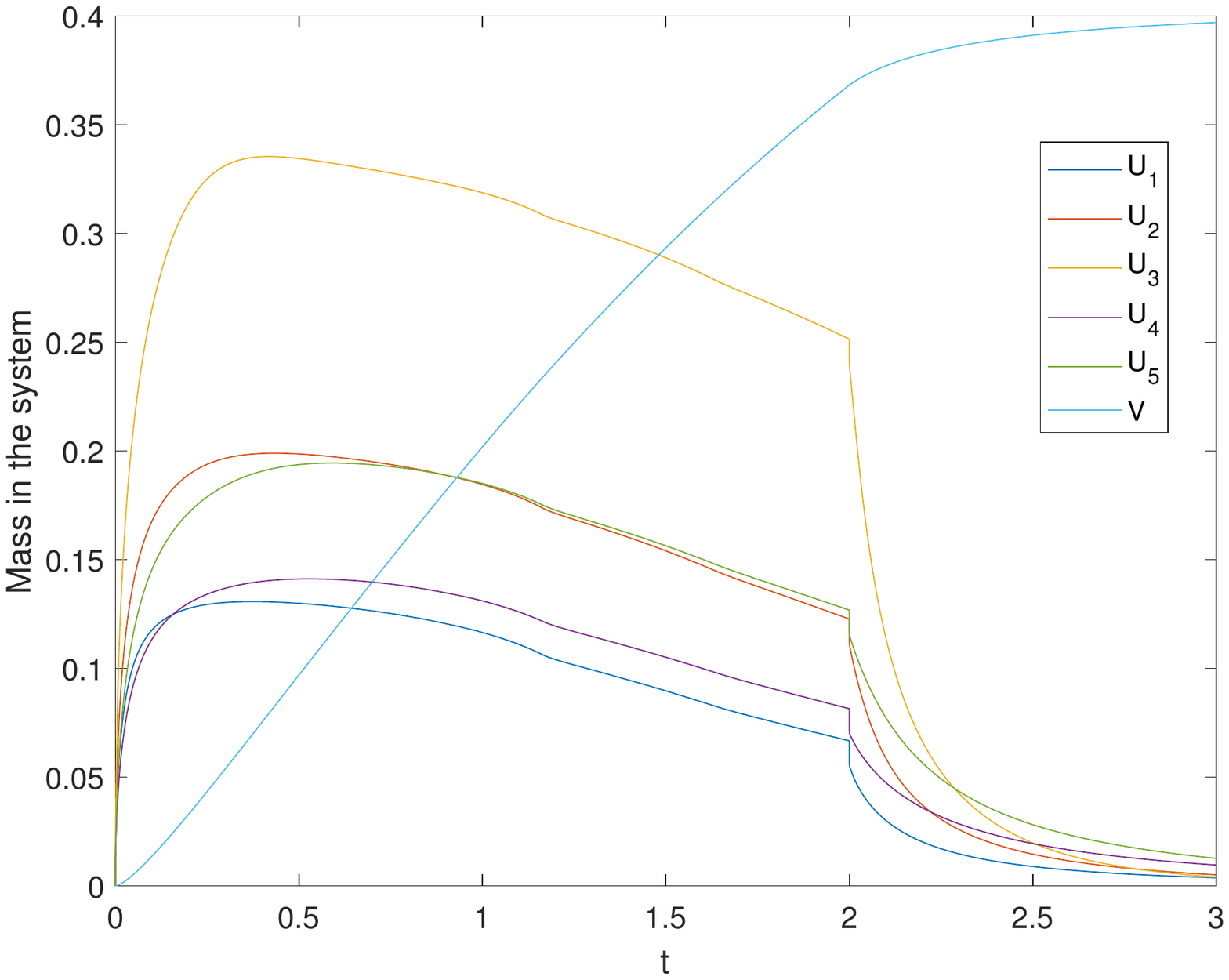}
  \caption{\it Mass of the system against time for the various components of the system in the case of 5 mobile species}
\label{Figsim5}
\end{figure}
\begin{figure}[h!]
\begin{center}
  \includegraphics[bb= 0 0 600 420, scale=.5]{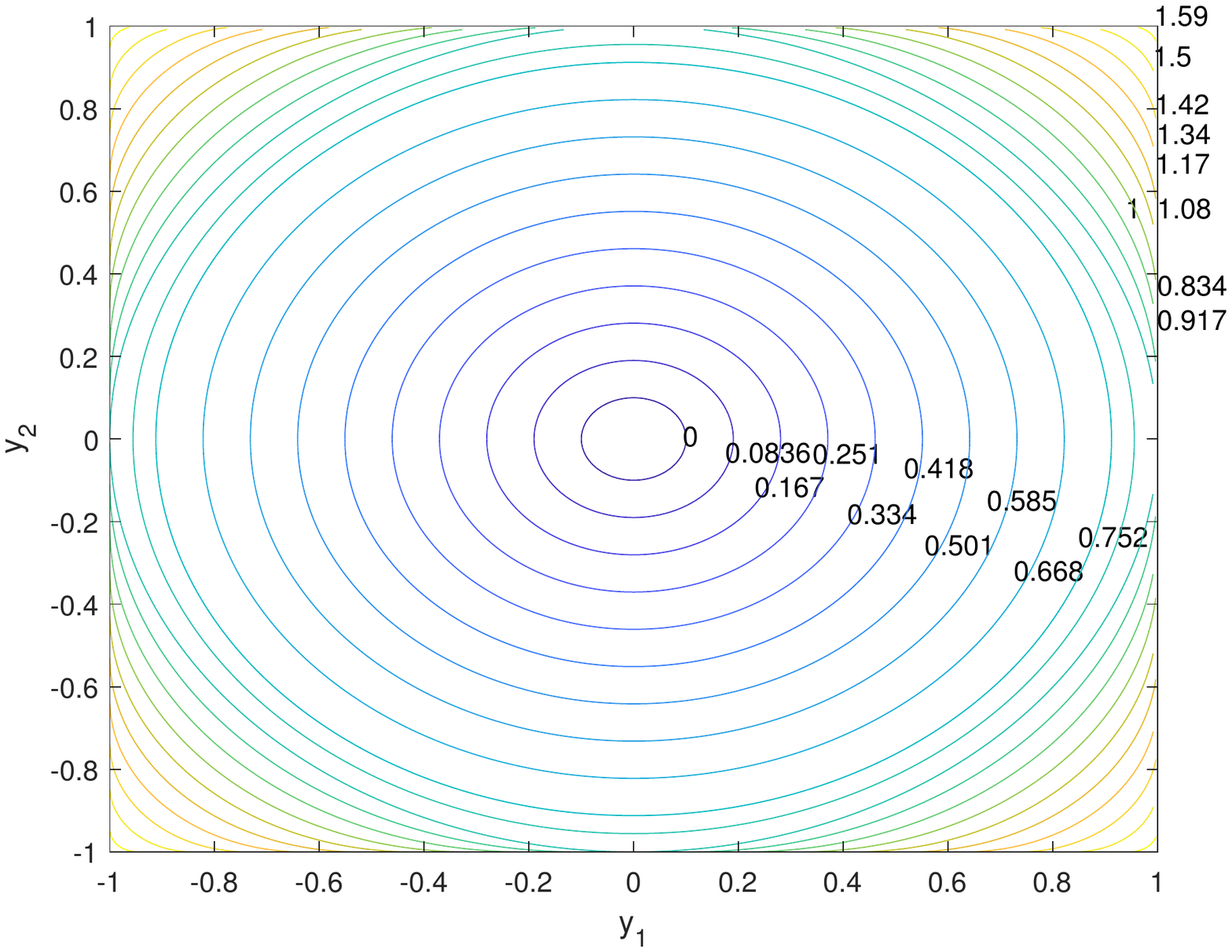}
\end{center}
\caption{\it Evolution of the moving boundary of the core formed by the immobile species for the case of configuration A at the point $x=0$. The numbers upon the curves indicate the time.}\label{cloggA}
\end{figure}
\begin{figure}[h!]
\begin{center}
  \includegraphics[width=11cm, height=6.5cm, bb= 0 0 550 220, scale=.7]{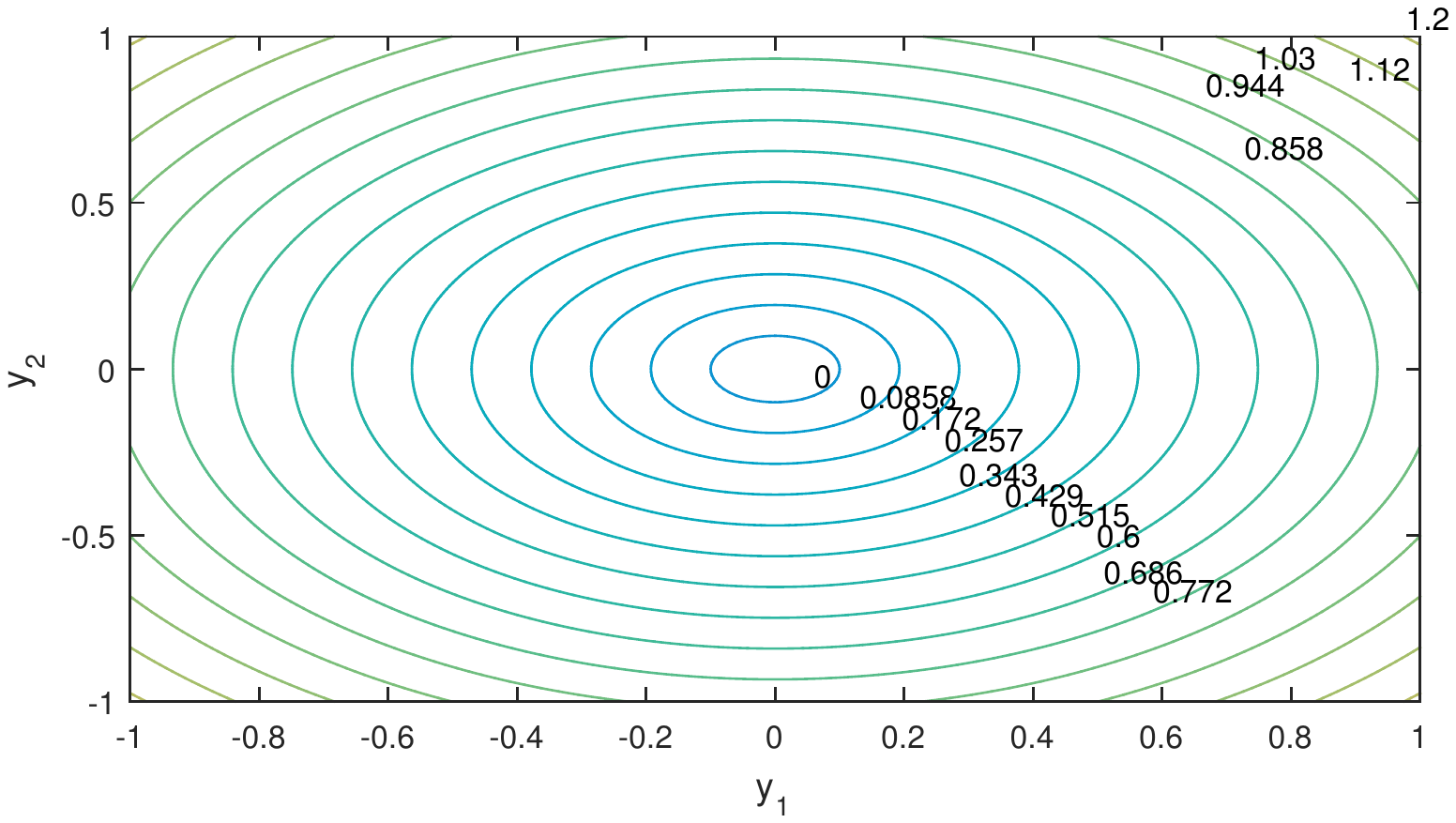}
\end{center}
\caption{\it Evolution of the moving boundary of the core formed by the immobile species for the case of configuration B. The numbers upon the curves indicate the time.}\label{cloggB}
\end{figure}

%____________________________________
In the same simulation we can also observe the evolution of the free boundary inside the cells at the point $x=0$.
More precisely, in the next set of figures, we  see the evolution of the core consisting of immobile species for a sequence of time steps.  In Figure \ref{cloggA} this effect is shown for the case that configuration A is assumed. As far as $r(x,t)<1$ the core has a cyclic  form and then for later times the corners of the cell are gradually filled by the immobile species with its boundary cyclic and tangential to the cell boundary.
%-----------------------------------------------------------------------------
%____________________________________
%\begin{figure}[h!]
%  \centering
%  \includegraphics[width=10cm ,height=10cm]{clogg5A.eps}
%  %width=1\textwidth
%  \caption{\it Evolution of the moving boundary of the core formed by the immobile species for the case of configuration A at the point $x=0$. The numbers upon the curves indicate the time.}
%\label{cloggA}
%\end{figure}
%____________________________________
%----------------------------------------------------------------------
Similarly for  the configuration B, the evolution of the moving boundary for various times is shown in Figure \ref{cloggB}.
%----------------------------------------------------------------------------------
%----------------------------------------------------------------------------------
 %_____________________________________________________________ 
%\begin{figure}[!htb]
%   \begin{minipage}{0.48\textwidth}
%     \centering
%%     \includegraphics[width=6cm,height=6cm]{clogg5A.eps}%\vspace{-6cm}
%     \caption{\it Evolution of the moving boundary of the core formed by the immobile species for the case of configuration A at the point $x=0$. The numbers upon the curves indicate the time.}\label{cloggA}
%   \end{minipage}%\hfill
%   \hspace{.5cm}
%   \begin{minipage}{0.48\textwidth}
%     \centering%\vspace{-3cm}
%  %  \includegraphics[width=6cm,height=6cm]{clogg5B.eps}%\vspace{-3cm}  
%          \caption{\it Evolution of the moving boundary of the core formed by the immobile species for the case of configuration B. The numbers upon the curves indicate the time.}\label{cloggB}
%   \end{minipage} %\caption{}
%\end{figure}
%________________________________________________________

%____________________________________
%\begin{figure}[h!]
%  \centering
%  \includegraphics[width=10cm ,height=10cm]{clogg5B.eps}
%  \caption{\it Evolution of the moving boundary of the core formed by the immobile species for the case of configuration B. The numbers upon the curves indicate the time.}
%\label{cloggB}
%\end{figure}
%____________________________________
%______________________________________________________
The numbers that are used  in these graphs, Figures \ref{cloggA} and \ref{cloggB}, are the same as those used in producing Figure \ref{Figsim1}.

Furthermore, in the next Figure \ref{varkappa} a simulation is presented for various values of the surface Thiele modulus $\kappa$.  Keeping the mass constant at the left boundary $x=0$ with increased   $\kappa$  ($\kappa=10$), we have more mass supply in the system and therefore higher values for $U_i$'s and $V$. The opposite effects is apparent for lower $\kappa$,  ($\kappa=0.1$).

\begin{figure}[h]
  \centering
    \includegraphics[bb= 0 0 500 320, scale=.6]{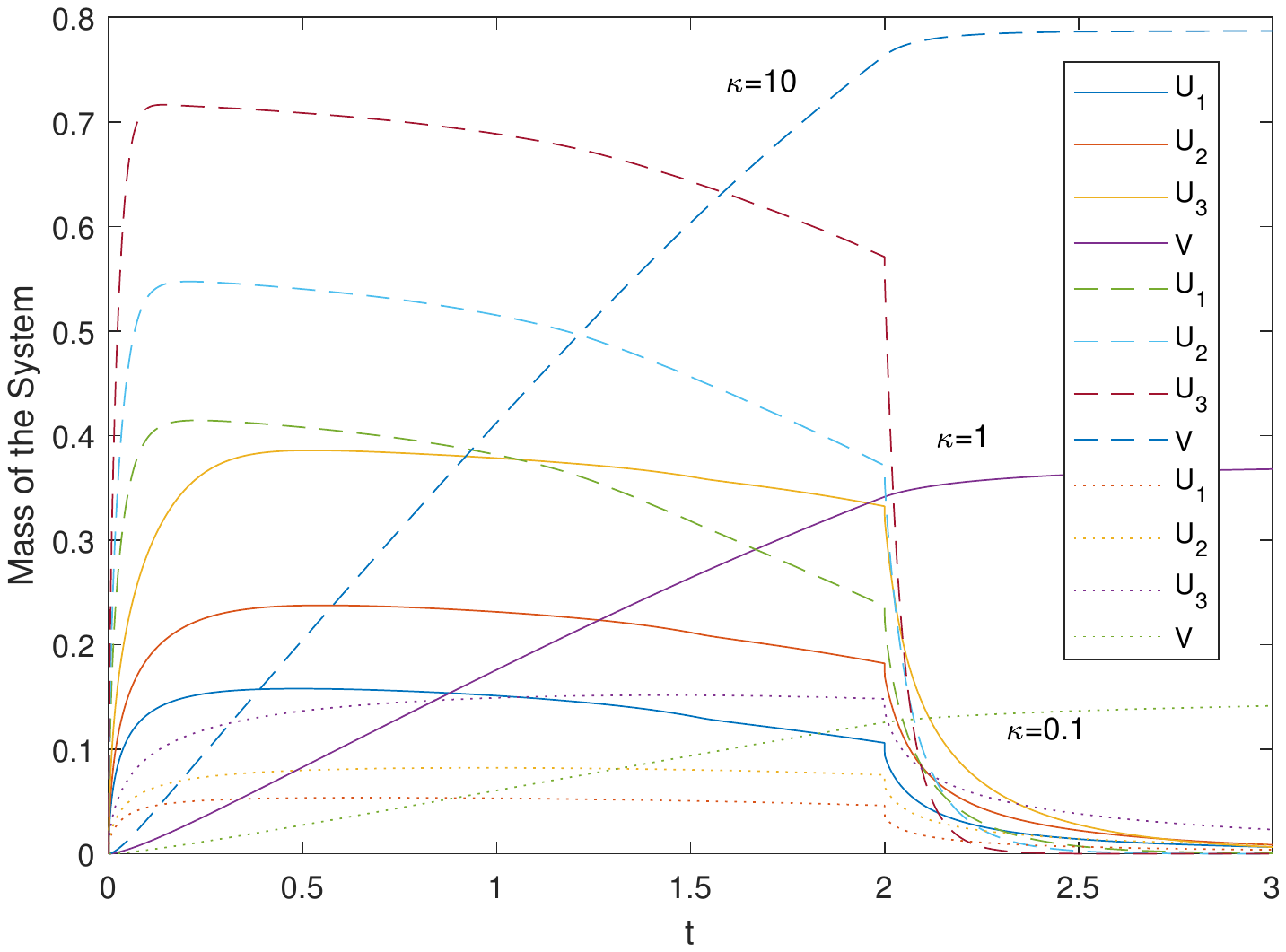}
  \caption{Plot of the system mass,  $U_i$'s and $V$, $i=1,2,3$,  for various values of $\kappa$, (a) $\kappa=10$ (dashed lines) ,  
  (a)  $\kappa=1$ (solid lines), (c) $\kappa=0.1$ (dotted lines). The rest of the parameters are the same as in Figure \ref{Figsim1}.
   }
  \label{varkappa}
\end{figure}

\begin{figure}[h]
  \centering  
    \includegraphics[bb= 0 0 600 320, scale=.5]{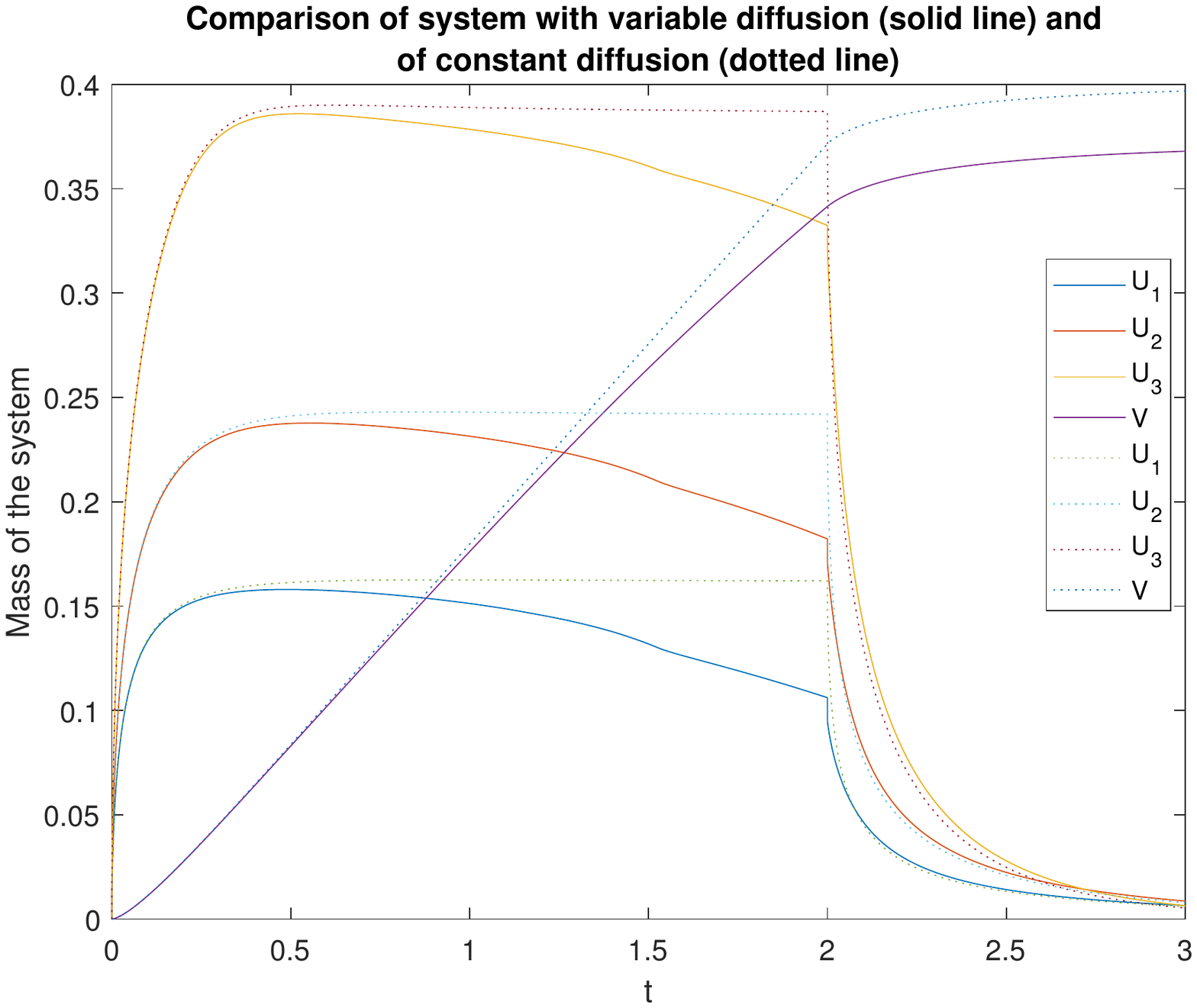}
  \caption{The breakthrough curve comparison for constant diffusion
    vs. pore-growth diffusion.  We can observe an initial agreement in both cases and then a decrease in the case of variable (due to pore growth) diffusion which is the result of the decreasing diffusion coefficient.
    However, after time $t=2$, where the boundary condition is dropped to zero the profile evens out to match the one of the larger constant diffusion.}
  \label{compDA}
\end{figure}

Finally, a comparison with the same model but with constant diffusion,  the same for all species, is shown in Figure \ref{compDA}.
Here we can observe that for the same species when we consider 
variable diffusion coefficient,  decreased with time, we have less amount of species diffusing inside the bulk of the material and therefore the total mass  of the  species is smaller compared with the case that we run the same simulation but with constant diffusion. When the supply from the left boundary stops at time $ t_0=2$, then the remaining mobile species collide to the core and finally the mass of the system in both cases evens out  (in the sense  the $u_i$'s and $V$ tend to a steady state). This behaviour seems to be robust.
 %

%________________________________HHH
%_______________________
%_____________________________________________________________ 
%\begin{figure}[!htb]
%   \begin{minipage}{0.48\textwidth}
%     \centering
%  %   \includegraphics[width=6cm,height=6cm]{compDAn.eps}%\vspace{-6cm}
%     \caption{\it The breakthrough curve comparison for constant diffusion
%    vs. pore-growth diffusion.  We can observe an initial agreement in both cases and then a decrease in the case of variable (due to pore growth) diffusion which is the result of the decreasing diffusion coefficient.
%    However, after time $t=2$, where the boundary condition is dropped to zero the profile evens out to match the one of the larger constant diffusion.}
%     \label{compDA}
%   \end{minipage}%\hfill
%   \hspace{.5cm}
%   \begin{minipage}{0.48\textwidth}
%     \centering%\vspace{-3cm}
%     \includegraphics[bb= 330 230 250 600, scale=.8]{Figporx.pdf}
%  %   \includegraphics[width=6cm,height=6cm]{Figporx.pdf}%\vspace{-3cm}  
%          \caption{\it Plot of porosity $\phi(x,t)$ against time $t$ at certain points $x_i$ in the domain $[0,1]$.\vspace{2.1cm}}
%          \label{Figporx}
%   \end{minipage} %\caption{}
%\end{figure}
%________________________________________________________

\begin{figure}[h]
  \centering
\includegraphics[bb= 0 0 500 320, scale=.6]{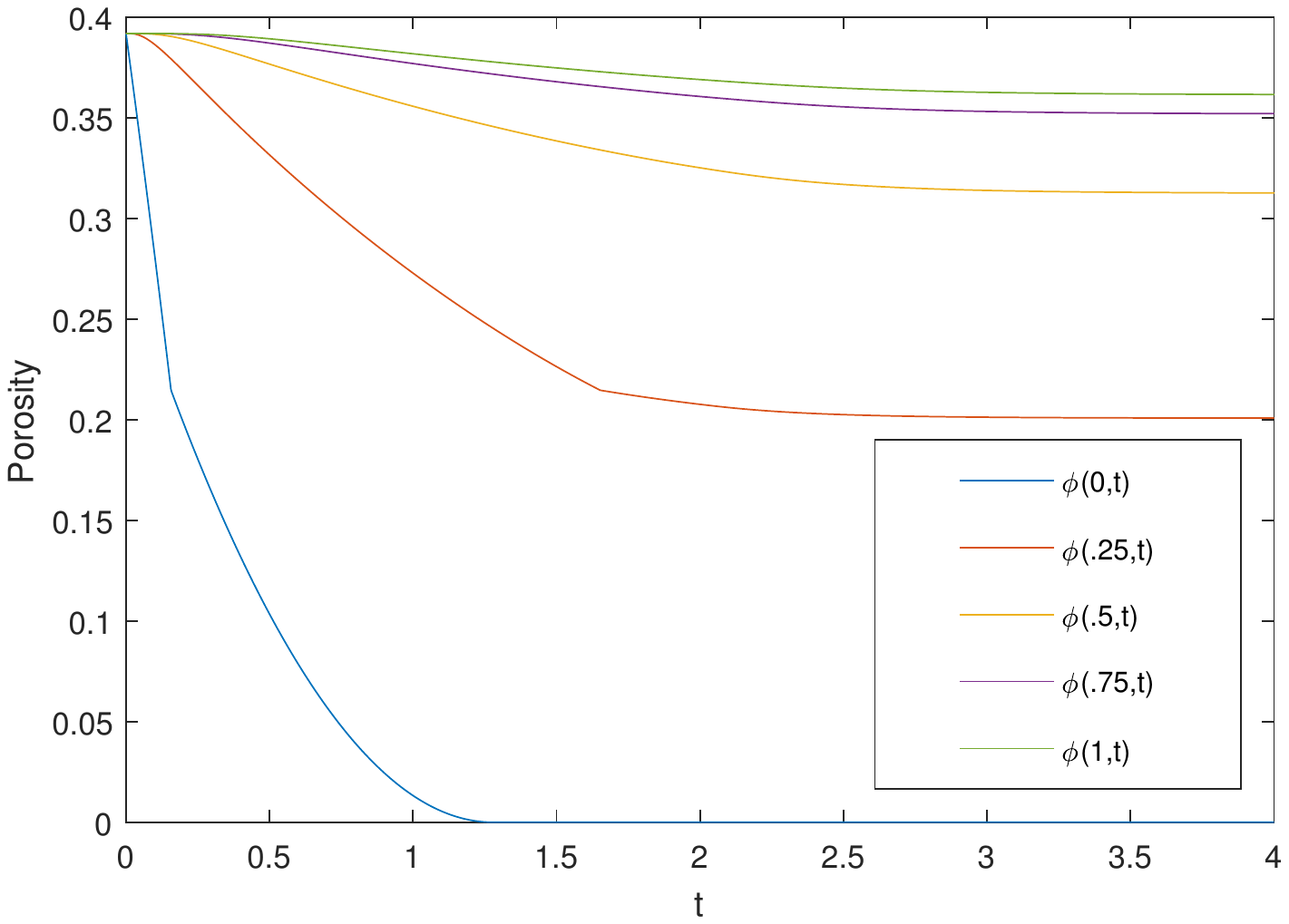}
  \caption{Plot of porosity $\phi(x,t)$ against time $t$ at certain points $x_i$ in the domain $[0,1]$.}
  \label{Figporx}
\end{figure}

Moreover, the porosity evolution is apparent in Figure \ref{Figporx}. The porosity function $\phi(x,t)$ is plotted here against time for certain values of the variable $x$ for $r(x,0)=.88$ ( $\phi(x,0)\simeq 0.39$) and with the rest of the values as in Figure \ref{Figsim1}. The porosity decay is much faster close to $x=0$
and actually becomes zero at time $t\simeq 1.22$ which means that we have clogging at that point and any boundary condition there does not effect any more the rest of the domain.
 Additionally, due to the latter fact, almost immediately after clogging the system settles down to  stationary state, close to that state observed at clogging time, since now we have the same problem  with no flux condition at both boundaries. 
%_______________________HHH
%_______________________

\begin{figure}[h]
  \centering
  \includegraphics[bb= 0 0 450 250, scale=.7]{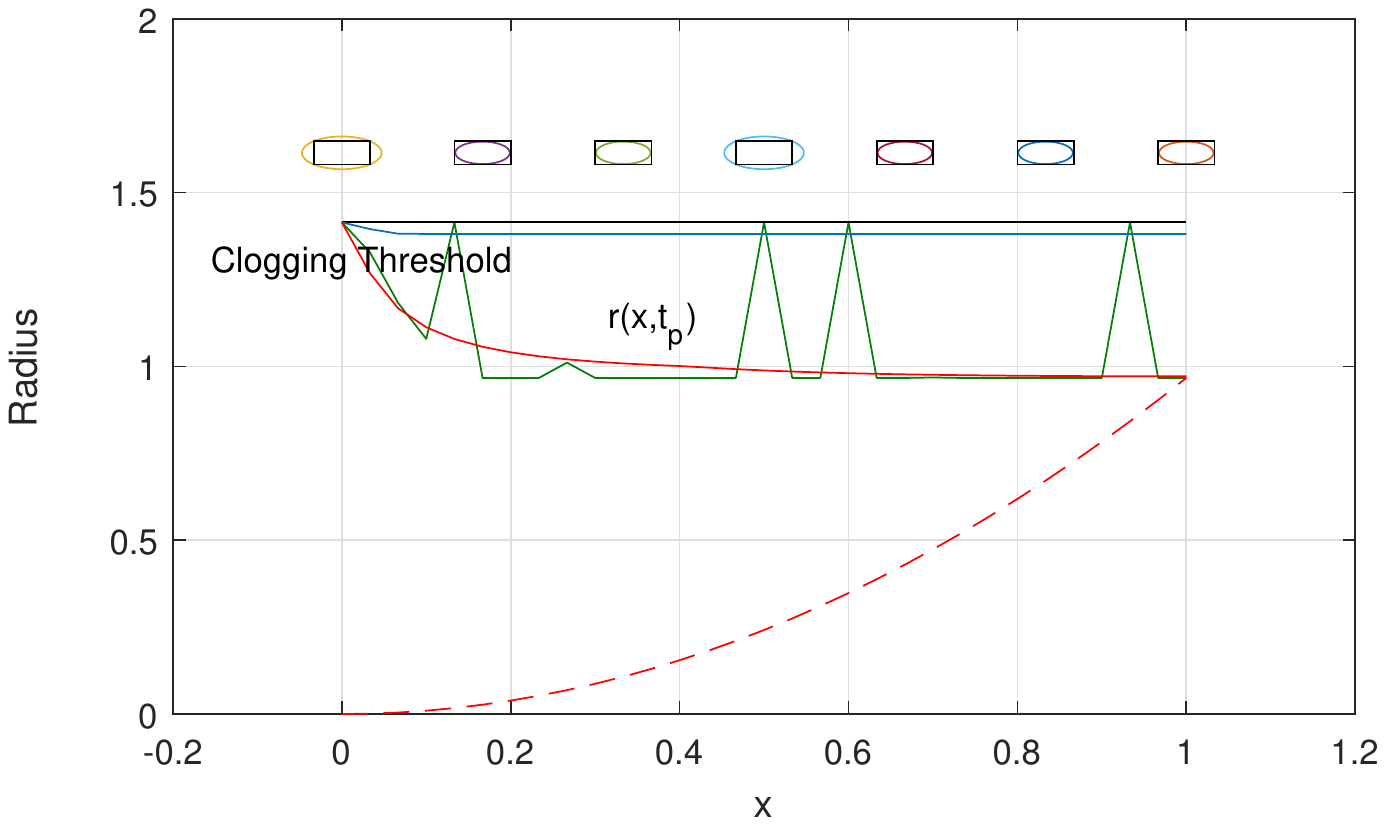}
  \caption{Form of radius distribution at time $t_p=2.5$ for the case that we take (i) $r(x,0)=x^2$ (red dashed line), (ii) $r(x,0)=1.38\,x^2$, (iii)  random normal distribution with mean $0.3$ and variance $0.8$. The radius $r(x,t_p)$ is plotted for these cases with (i) red line, (ii) blue line, (iii) green line. In the latter case at the top of the graph  a schematic image of some square cell is plotted together with corresponding circle of radius $r(x_i,t_p)$.}
  \label{fig_radius1}
\end{figure}

%_______________________
In the following set of simulations shown in Figure \ref{fig_radius1}, we investigate the effect of non-uniform  initial pore radii distribution, i.e. $r(x,0)=r_a(x)$ (non constant).
For the first simulations, we take $r(0,t)=x^2$ and for the case when configuration A is assumed (red dashed line). The radius $r(x,t)$ is plotted for time $t=t_p=2.5$ for $t_p=\frac56 T$ ($[0,T]$ is the simulation interval with $T=3$) against $x$ (red solid line).
 Note that the time $t_p$ is chosen appropriately at each of the following simulations so that to capture the clogging behaviour of the system
We have clogging when $r$ reaches $\sqrt{2}$ at some point and this threshold is denoted by the black line. When $r$ reaches $\sqrt{2}$ at some point $x_0$ then the pores in that point are clogged and we have no mobility from the area $x<x_0$ to the area $x>x_0$. For this case this happens initially at the point $x=0$.

Changing the initial distribution of $r$ by taking now $r(x,0)=r_0 x^2$ for $r_0=1.38$ close to $\sqrt{2}$, we may have a situation that the whole domain is almost clogged simultaneously ( blue line in Figure \ref{fig_radius1}).

Finally, we obtain an  interesting situation  if we take the values of $r(x,0)$  at the points of the partition given by a random distribution. In this case starting with a normal distribution with mean $0.3$ and variance $0.8$  we obtain at time $t_p$ a profile (green line) in which clogging is exhibited at the points $x\simeq 0,\,0.1,\, 0.5,\, 0,6,\, 0,9$. Also for this case and this specific time step $t_p$, at certain selected points $x_i$, a schematic image of a square cell is plotted together with circle of radius $r(x_i,t_p)$.
Starting with a  continuous uniform (rectangular) distribution in the interval $[0,\sqrt{2}]$ we obtain a similar result.

The occurrence of clogging can be seen in these graphs when $r(x,t_p)$ reaches 
$\sqrt{2}$ at some point $x$ and time $t_p$. Now, one could also ask: {\em How can one distinguish between through-pores and ink-bottled pores}? Our current level of understanding is that 
if $r$ is lower than 1, then we have through-pores, while for $r>1$ we have bottle-through pores.  Looking at Figure  \ref{fig_radius1}, we see that (at a fixed time slice) the clogging threshold is reached in a number of points (e.g.  $x=0.5$ and $x=0.6$). In this region we expect blind and closed pores  to appear\footnote{ From a catalysis-oriented perspective, e.g. cf. \cite{Leofantia}, pores can be closed (not accessible from the air parts of the pore matrix), blind (open at only one end), or through (open at both ends). Each pore can be isolated or, more frequently, connected to other pores to form the porous network. The irregular shape of the pores and their connectivity cause a molecule or a colloidal particle to cover a distance greater than the grain size when passing through the porous material.}.

\begin{figure}[h]
  \centering
  \includegraphics[bb= 0 0 380 200, scale=.7]{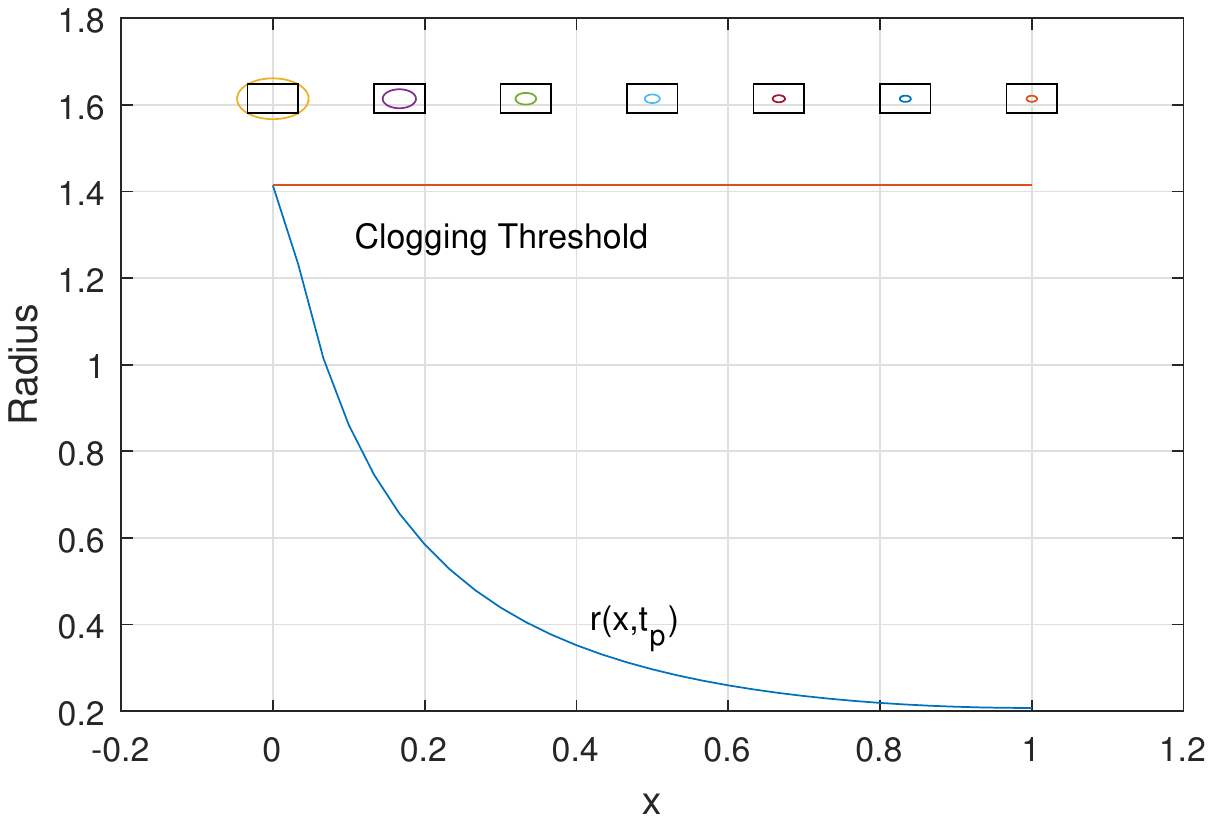}
 \caption{Form of radius distribution $r(x,t_p)$ at time $t_p=2$ for the case of  configuration B and having initially $r(x,0)=x$. We observe clogging at the point $x=0$. At the top of the graph, a schematic image of square cells centered at the points $x_i=(i-1) h$, $h=0.142$, $i=1\ldots 7$,
   is plotted together with the corresponding circle 
of radius $r(x_i,t_p)$.} 
  \label{fig_radius2}
\end{figure}

A similar experiment is demonstrated in Figure \ref{fig_radius2} but for the case of configuration B. The initial distribution here is taken to be linear $r(x,0)= x$ and here $t_p=2$. In the top of the figure we can see the state of the cells at this time step. Since material is inserting from point $x=0$, clogging takes place initially  there and then we see no significant variation  in the evolution of the process at later times.

    Based on the numerical results shown in Figure \ref{Figporx}, Figure \ref{fig_radius1} and Figure \ref{fig_radius2}, we conjecture the following: 
    \newline
    \fbox{
    \parbox{0.97\textwidth}{
    If ink-bottle pores are clogged, then  
     the storage capacity of the material reaches a certain saturation threshold simultaneously with a potential  
     increase of the averaged transport up to a maximum level. On the other hand, if through pores become clogged, then the averaged transport is decreased leading to a raise in the storage capacity. 
    }}

      In Section \ref{sc}, we bring in some further partial support for this conjecture, by estimating numerically the mass of colloids that have deposited on the pores surfaces and were not depleted from there during a given time interval.

    \subsection{A posteriori estimation of the storage capacity}\label{sc}

     A major question related to modern porous media applications is how much ${\rm CO_2}$,  energy, colloids,  residual nuclear particles etc. can be stored in the fabrics of a given heterogeneous medium; see \cite{Lou} for instance.  If, like in our case, the distribution, shape and volume of the microstructures is controlled, then the storage capacity of the porous medium can be estimated numerically employing the upscaled model equations.  Of course, even in our rather simple case,  such a task is quite complex  especially if one targets accurate estimations. The success depends 
     on a lot of factors such as the  space inside the pores available for the colloids to accumulate, the void space inside the pores where colloids feel free to move, number and location of clogging regions, etc. Part of this information is incorporated in our model  in the functions $L(x,t)$  (tortuosity information,  in the sense that with increasing $L$ the average travelling length of a particle 
     in the non-solid part of a cell also increases) and the porosity $\phi(x,t)$. In addition to these aspects, additional  factors that are not included in this model such as the shape and size of the colloids, the way that they are packed and accumulated in the pore surface etc. should be taken into account if one wants a better insight. On top of this, due to the nonlinearity and coupling in the upscaled system, we cannot aim to determine  a specific formula of the storage capacity of the system.  However, in our context we can give an estimation in an indirect way.  If we focus at a point of our domain at time $t=0$ we may  consider that it has minimum storage capacity. As the system evolves with time, the pores associated to this macroscopic spatial point tend to fill out, so more matter (colloids) can be stored in the system. We refer to this scenario as an increase in the  storage capacity of the medium until the time where clogging occurs. Then the capacity reaches a threshold.
    The latter means that we have no void space or free surface at this point i.e. the pores are either vanished or closed and hence additional storage becomes impossible. Based on such picture, we can define a local capacity storage indicator as well as a global capacity storage indicator  by exploiting the difference $SC_\ell(x_c,t)= (v(x_c,t)-v(x_c,0)$, where $x_c$ is the point where clogging occurs. 
    % indicative function of the storage capacity at each time and at a point $x_c$ should be the (final) amount of $v$ when clogging occurs minus the amount the accumulated material at time $t$, i.e. $SC(x_c,t)=v(x_c,T)-v(x_c,t)$. 

    Also, for the case that at time $t=0$ we have a uniform initial distribution of pore radii
    and since the geometry in the microstructure is considered uniform,
    for any another point say $x$ (even when not clogging occurs there) we expect that we should have the same maximum storage capacity $v(x_c,T)$ and its (local) storage capacity at time $t$ should be $SC_\ell(x,t)=v(x,t)-v(x,0)$. A global indicator of the capacity storage is then $SC_g(t)= \int_\Omega (v(x,t)-v(x,0))dx$.

\begin{figure}[h!]
      \centering
       \includegraphics[bb= 0 0 450 280, scale=.7]{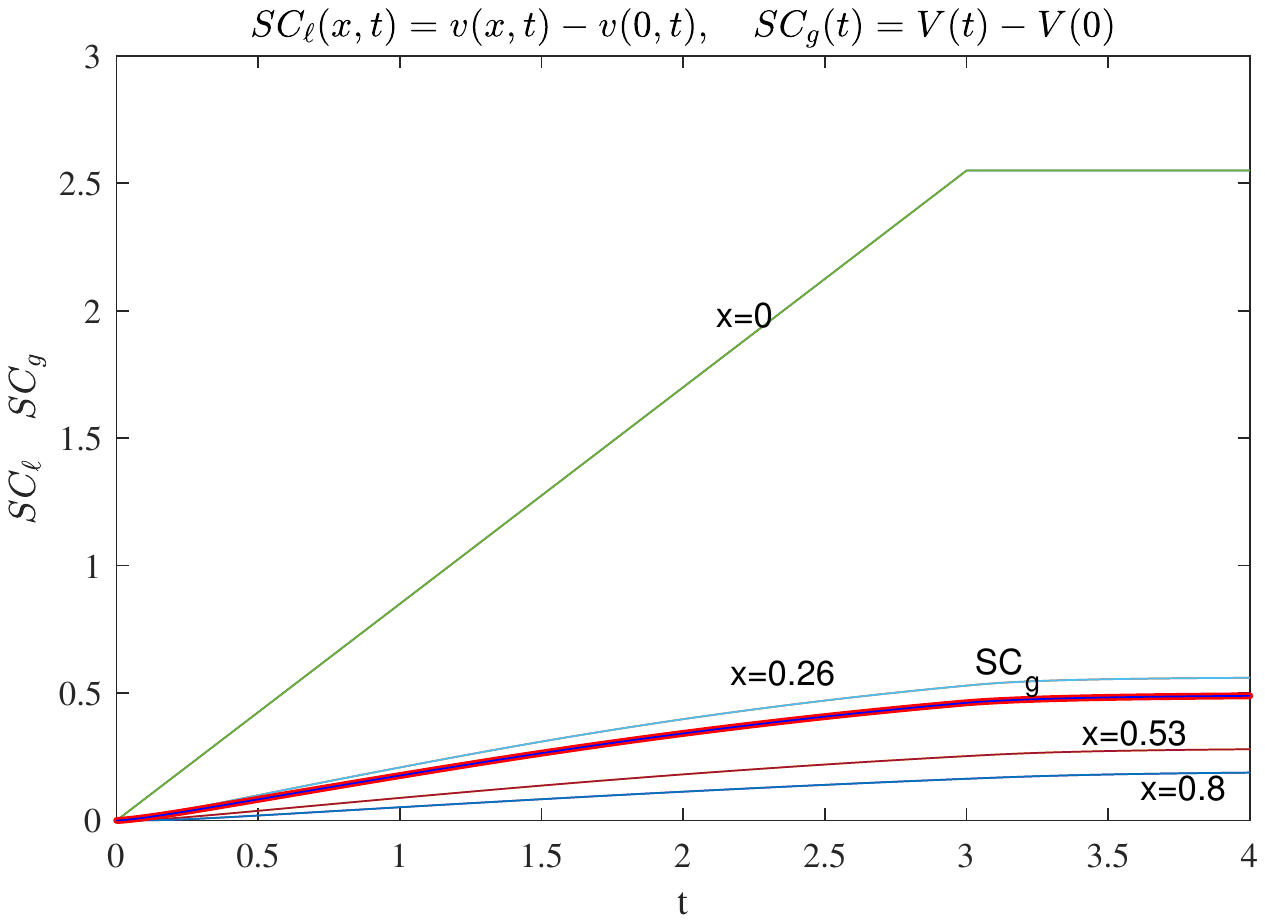}
     \caption{Estimation of the local and global storage capacity  indicators as functions of time, where $V(t)=\int_\Omega v(x,t)dx$.} 
      \label{fig_stcap1Abn}
    \end{figure}

\begin{figure}[h]
      \centering
       \includegraphics[bb= 0 0 580 300, scale=.5]{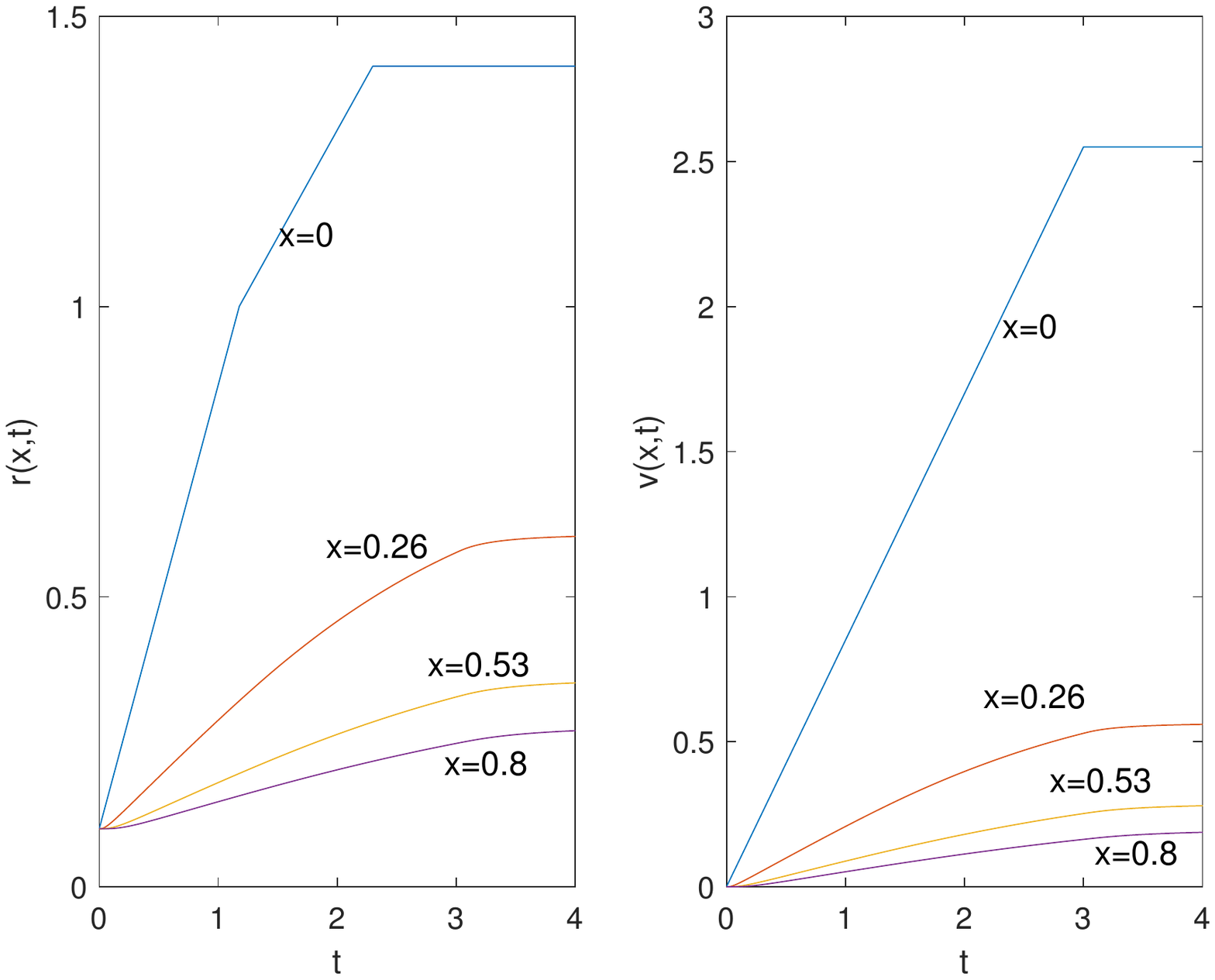}
     \caption{The pore radius $r(x,t)$ and the concentration profile $v(x,t)$ are plotted  with respect to time for $x=0, \,0.26,\,0.53,\,0.8$.} 
      \label{fig_stcap2A}
    \end{figure}

\begin{figure}[h!]
   \begin{center}
      \includegraphics[bb= 0 0 450 350, scale=.6]{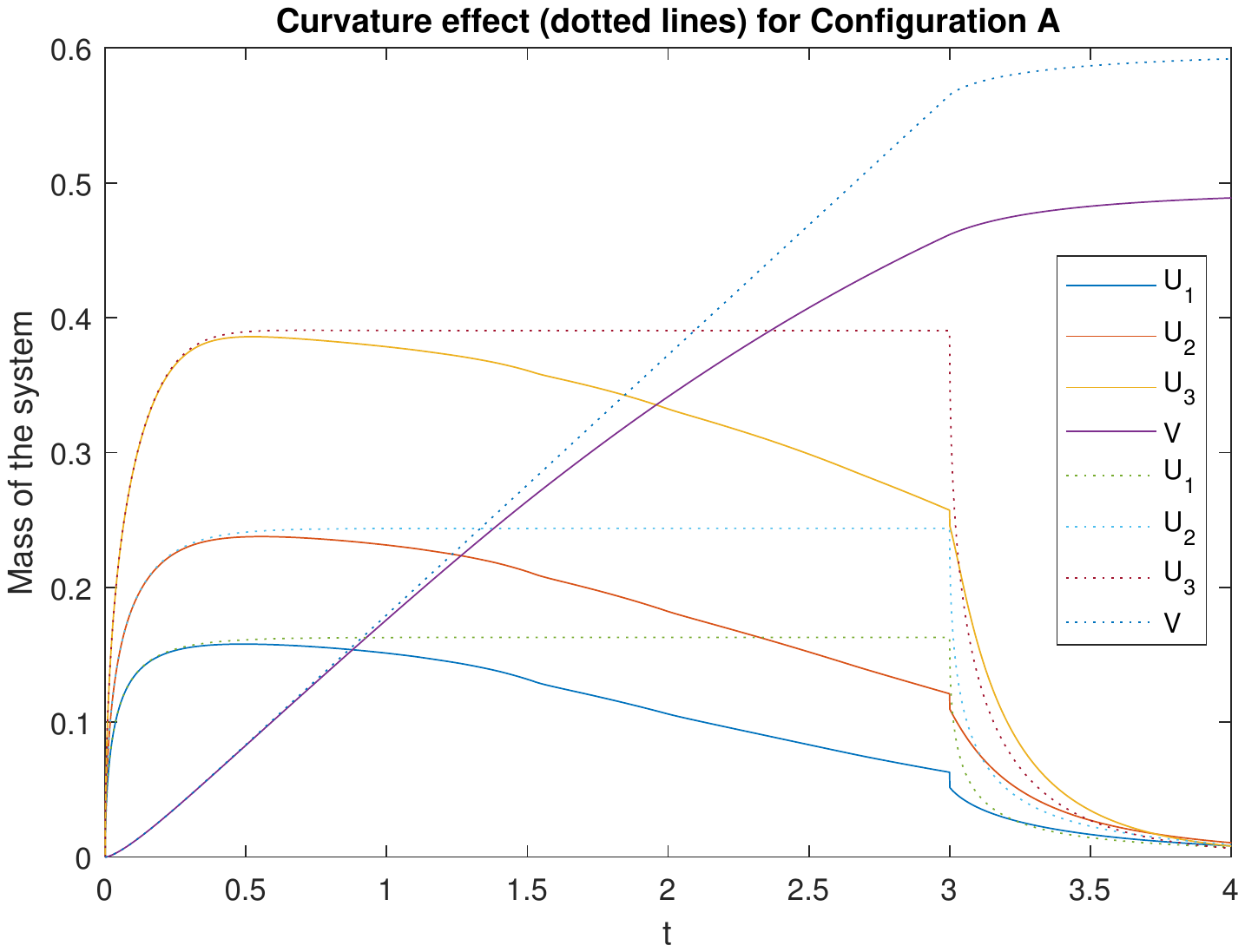}
\end{center}   
     \caption{Total mass in the system for configuration A, when the term $-\alpha\frac{1}{r(x,t)}$, $\alpha=0.1$  is added to the equation of $r'(x,t)$ as given by \eqref{eq:model-dim-4}.} 
      \label{fig_curvature}
    \end{figure}

    To get further insight, in the numerical experiment shown in Figure   
    %\ref{Figsim1} 
    \ref{fig_radius2},  we extend the inflow at the boundary $x=0$ for timeslot $[0,t_0]$, $t_0=3$ (cf. (\ref{mod1b})) 
     so that enough mass is inserted in the system allowing us  to 
    observe clogging at the point $x=0$. In  Figure \ref{fig_stcap1Abn}, the local storage capacity indicator 
    $SC_\ell(x,t)$ is plotted vs. time for various choices of spatial points $x$. The thick line in the picture represents the global storage capacity indicator  vs. time.
    In the next plot (see Figure \ref{fig_stcap2A}) the pore radius $r(x,t)$ (left) and the concentration of the immobile $v(x,t)$ (right) species are plotted against time for the same choice of spatial points. 
   The similarity between the two behaviors points out the fact that we expect that any eventual closed form formula for a storage indicator should incorporate a suitable monotonic dependence on the pore radii. More investigations are needed to shed light on this effect.
      One could also wonder what would be the effect of an additional curvature term like $-\alpha\frac{1}{r(x,t)}$ ($\alpha>0$) in the definition of the speed of the microscopic free boundaries.  Interestingly, our simulations show that, for both configurations A and B, the storage capacity of the porous medium increases in the sense that  more mass is stored in the system compared with the no-curvature case; see Figure \ref{fig_curvature} for an illustration of what happens in the case of configuration A; for configuration B a similar result can be pointed out. To understand why and how this effect occurs, one needs to study the time evolution of more types of shapes, preferably less regular and non-symmetrical, but this goes beyond our aim here.
 
%\newpage
\section{Conclusions} 

A few conclusions can be drawn for this framework:
%\begin{itemize}

If the internal structure of a porous material is fixed with respect to time, then the effective transport and storage properties  of the material can usually be sufficiently well estimated by empirical measurements.   Efficient mathematical averaging methodologies are needed  to deliver realistic quantitative estimations of the transport and storage capacity of heterogeneous materials with evolving internal structures (think of the effect of freezing or thawing ice lenses in cementitious-based materials) - for such situations, experimental techniques fail to be applicable.

 What concerns the choices of microstructures considered here,  we have seen that both the volumetric porosity and the surface porosity 
(i.e. the ratio of the area of the voids in a plane cross section of the porous medium to the total area of the cross section),
are  proportional to $A(x,t)$,
  the effective reaction constant (via Henry-like exchange) is  proportional to  $L(x,t)$, and the tortuosity as well as the effective  diffusion coefficient to $\left( 1+ \frac{1}{A(x,t)}\int_{Y(x,t)} \frac{\partial w_1}{\partial y_1}dy \right)$.
  For the simplified setting  taken to describe the large radii case, $w_1$ is a  cell function which together with the analytic expressions for $A(x,t)$ and $L(x,t)$, depends on the selected microscopic configuration ($A$, $B$, or some other admissible situation). 

 We do not have a way to express analytically the effective storage capacity of a heterogeneous domain  (with or without evolving microstructures!). However, for the particular choices of microstructures chosen here we are able to estimate numerically an averaged transport flux allowing for local reductions eventually due to clogging or to a sudden formation of localized high-density zones induced by the microstructure evolution.   Using the evolution equation for  the mass of deposited colloids, we can define and compute numerically local and global capacity storage indicators. 

If a more complex shape of locally-periodic microstructure needs to be approximated, then a front-tracking approximation scheme (e.g. using a suitable level-set equation)  will need to substitute our numerical approach. The same holds if $\Omega$ is two or three dimensional.       Should a 3D setting be treated, then (\ref{eq:model-dim-4}) has to be changed accordingly. We expect that in this case the estimation of the storage capacity is more difficult to obtain especially if ball-like grains are replaced by other growing shapes.

 In spite of the fact the the FBP is well-posed locally in time (cf. Appendix \ref{app}) and that we have reasons to expect that the macroscopic model is globally well-posed (e.g. trusting the approach from \cite{ray2013drug,ray2015}), we are unable at this moment to perform  the rigorous homogenization for this setting. The missing ingredient is a lower uniform in $\varepsilon$ estimate of the local existence time of the oscillating FBP.  The weak solvability of the microscopic system can be handled up to the first  clogging point. Looking at the macroscopic equations, if clogging occurs transport (elliptic part) degenerates simultaneously with other terms in the upscaled equations. A similar situation has been encountered for clogging pipe-like pores in \cite{Noorden}. In both these scenarios, one lacks a rigorous upscaling of the degenerate evolution equations. 
 
  In this framework, we have discussed only the case when convective flows are absent. Further investigations are needed to clarify how flows with moderate oscillatory speeds can be handled in this setting. 
 
%\end{itemize}

%___________________________________________________________________________

\section{Acknowledgments} 

The authors thank Dr. Oleh Krehel (Eindhoven) for helpful discussions.  AM acknowledges  partial financial support from NWO-MPE ``Theoretical estimates of heat losses in geothermal wells" (grant nr. 657.014.004)   and from Swedish Research Council ``Homogenization and dimension reduction of thin heterogeneous layers" (grant nr. 2018-03648).

\appendix

\section{Comments on the local weak solvability of the original FBP}\label{app}

We consider the simplified case when the solid phase in $\Omega$ consists of a collection of $M\in \mathbb{N}$ disjoint balls with centers in $x_\ell\in \Omega$ ($\ell\in\{1,\dots,M\}$) placed on a regular $\delta$-grid (see Figure \ref{fig:locally-periodic});  the grid's cells have diameter $\delta>0$.  We look into the case when $0<\epsilon\ll \delta$ and assume that the centers do not move and the ball-shape is preserved during the evolution of the free boundary. We fix the notation as follows:
$$ \Omega_0^\varepsilon(r_\varepsilon(t)) := \cup_{\ell=1}^M \varepsilon Y_0(x_\ell, r_\ell(t)) \mbox{ and }   \Gamma^\varepsilon (r_\varepsilon(t)):=\cup_{\ell=1}^M \epsilon \partial Y_0(x_\ell, r_\ell(t)).$$ 
 To simplify the writing, let us refer  to $\Omega_0^\varepsilon(r_\varepsilon(t))$ and to  $\Gamma^\varepsilon (r_\varepsilon(t))$ simply as  $\Omega_0^\varepsilon(t)$  and respectively, $\Gamma^\varepsilon (t)$.
Consequently, $\Omega^\varepsilon(t):=\Omega-{\rm clos}(\Omega_0^\varepsilon(t))$,  $\Omega_T^\varepsilon:=\cup_{t\in (0,T)}\left(\Omega_0^\varepsilon(t)\times \{t\}\right),$ and $\Omega_T:=\Omega\times (0,T)$ for some $T>0$.

Further, let the maximum existence time interval be denoted by
$$t^{\delta\varepsilon}:={\rm inf}\{t: 2\varepsilon r_\ell(t)=\delta, \ \ell\in\{1,\dots, M\}\}. $$
Note that, by construction, we have $t^{\delta\varepsilon}>0$.  Set $u^\varepsilon=(u_1^\varepsilon,\dots, u_N^\varepsilon)$ for the vector of colloidal populations  and $r^\varepsilon(t)=(r_1^\varepsilon(t),\dots, r_M^\varepsilon(t))$ for the vector of radii of the growing balls.   

The FBP we are considering is: Find  the triple $(u_\epsilon,v_\epsilon,r_\epsilon)$ satisfying the following evolution system:
\begin{align}
  \partial _t u_i^{\varepsilon}&=\kappa d_i \Delta u_i^{\varepsilon} + R_i(u^{\varepsilon})      &&\text{in }\Omega_0^\varepsilon(t) ,i\inRange{N},\\
  -d_i\nabla u_i^{\varepsilon}\cdot \normal_\varepsilon &=\varepsilon 
  (\dca_i u_i^{\varepsilon}-\beta_i v^{\varepsilon}) &&\text{on } \Gamma^{\varepsilon}(t),\\ 
  \partial _t v^{\varepsilon}&= \eta\left(\sum_{i=1}^N \dca_i u_i^{\varepsilon}- \beta  v^{\varepsilon}\right) &&\text{on }\Gamma^{\varepsilon}(t),\label{eta1}\\ 
  & r_\ell^{\varepsilon}(t) {r_\ell^{\varepsilon}(t)}' = \alpha \int_{ \partial Y_0(x_\ell, \varepsilon r_\ell(t)) }\partial _t v^{\varepsilon} &&\text{in }\Omega_0^\varepsilon(t), \label{eta2}
\end{align}%\dcb
endowed with the following boundary and initial conditions
\begin{eqnarray}
\nabla u_i^{\varepsilon}\cdot \normal_\varepsilon =0 \quad \text{on } \partial \Omega_T,
\end{eqnarray}
\begin{align}
  u_i^{\varepsilon}(0,x)={u_i}^0(x) &&\text{in }\Omega_0^\varepsilon(t) ,i\inRange{N},  \\
   v^{\varepsilon}(0,x)={v}^0(x) &&\text{on }\Gamma^{\varepsilon}(t), \\
 r_\ell^{\varepsilon}(0)={r}_\ell^0, \ell\inRange{M}.
\end{align}
In (\ref{eta1}), the function $\eta:\mathbb{R}\to [0,\infty)$ defined via $\eta(s)=s^+$ is taken to ensure jointly with (\ref{eta2}) that the balls $Y_0(x_\ell, \varepsilon r_\ell(t))$ can only increase their size in time. The positivity of   $\partial _t v^{\varepsilon}$ is an essential aspect; it needs either to be imposed by employing $\eta$, or it must be proven to hold in a certain parameter regime (using sub- and super-solutions). %Here we only wish to illustrate a proof idea for local-in-time well-posedness of the FBP relying on a fixed-point argument of Banach type. 
It is worth mentioning at this point that the rigorous averaging of this setting is in progress. The main challenge  is at least twofold: (i) a suitable scaling of the balls as well as a proper dependence $\delta=\delta(\varepsilon)$ must be identified such that eventual boundary layers arising around the growing surfaces affect the averaging in a controlled way, (ii) the size of the maximum time interval where the local
 solvability of the FBP can be ensured is  depending in an unsuitable way on the choice of $\varepsilon$. Especially what concerns (ii), the situation here is totally different compared to the better understood case of Ostwald ripening. 
  
Here, we only sketch the proof idea for the application of the Banach-fixed-point theorem to this setting ensuring  local-in-time well-posedness of the FBP. We transform   the moving domain $\Omega_0^\varepsilon(t)$ into the fixed domain $\Omega_0^\varepsilon(0)$ by some diffeomorphism depending not only on  the choice of $r^\varepsilon(t)$ but also on its regularity and uniform boundedness in suitable norms; see also  
\cite{Karali12}, where the authors have fixed the free boundaries in a multiple-connected domain related to Ostwald ripening with kinetic undercooling. We refer the reader to \cite{Malte} where a simple connected domain has been treated. 

Take 
$$\Phi(\cdot, r^\varepsilon(t)): \Omega_0^\varepsilon(0)\to \Omega_0^\varepsilon(t) \mbox{  given by } \Psi(x,t)=\Phi(x,r^\varepsilon(t)).$$
Herewith, the newly transformed functions are $\tilde u^\varepsilon(x,t)=u^\varepsilon(\Psi(x,t),t)$ and  $\tilde v^\varepsilon(x,t)=v^\varepsilon(\Psi(x,t),t)$.   We now make use of the following notation referring to the (active part of the) structure of our heterogeneous medium in the fixed-domain formulation, viz. $\Omega_T^\varepsilon(0):=\cup_{t\in (0,T)}\left(\Omega_0^\varepsilon(0)\times \{t\}\right)$.

Given sufficiently regular $r^\varepsilon$ and $\frac{1}{r^\varepsilon}$, it results that the pair of functions $(\tilde u^\varepsilon,\tilde v^\varepsilon)$ must satisfy the parabolic problem

\begin{eqnarray}
&\sqrt{{\rm}det(D\Psi^TD\Psi)}\partial_t \tilde u_i^\varepsilon- k d_i {\rm div} \left(\sqrt{{\rm}det(D\Psi^TD\Psi)}(D\Psi^T D\Psi)^{-1}\nabla \tilde u_i^\varepsilon \right)= R_i(\tilde u^\varepsilon)+\nonumber\\
&+\sqrt{{\rm}det(D\Psi^TD\Psi)} D\Psi^{-T}\nabla \tilde u_i^\varepsilon\partial_t \Psi\mbox{ in }\Omega_T^\varepsilon (0) \label{ai}
\end{eqnarray}
\begin{eqnarray}
\nabla \tilde u_i^\varepsilon\cdot n =0 \mbox{ at }\partial\Omega_T^\varepsilon (0),
\end{eqnarray}
\begin{eqnarray}
D\Psi^{-T}\nabla \tilde u_i^\varepsilon = \varepsilon (a_i  \tilde u_i^\varepsilon -\beta_i \tilde v_i^\varepsilon)  \mbox{ at }\partial\Omega_T^\varepsilon (0),
\end{eqnarray}
\begin{eqnarray}
\sqrt{{\rm}det(D\Psi^TD\Psi)} \partial_t \tilde v^\varepsilon &=& \sum_{i=1}^N (a_i  \tilde u_i^\varepsilon -\beta_i \tilde v_i^\varepsilon) +\nonumber \\
&+& \sqrt{{\rm}det(D\Psi^TD\Psi)} D\Psi^{-T}\nabla \tilde u_i^\varepsilon\partial_t \Psi \mbox{ at }\partial\Omega_T^\varepsilon (0). \label{af}
\end{eqnarray}
Reasoning as in  Theorem 4.1 in Ref. \cite{Karali12}, standard theory of parabolic regularity ensures that if $(\tilde u_i^\varepsilon(0),\tilde v^\varepsilon(0))\in H^1(\Omega_0^\varepsilon (0))^N\times H^1(\partial\Omega_0^\varepsilon (0))$, 
 $r^\varepsilon\in W^{1,\infty}_+(0,T)^M$, and $\frac{1}{r^\varepsilon_\ell}\in L^\infty(0,T)$ for all $\ell\inRange{M}$ and some $T\in(0,t^{\delta\varepsilon})$, then problem (\ref{ai})--(\ref{af}) admits a unique solution 
$$(\tilde u^\varepsilon,\tilde v^\varepsilon) \in \left(L^\infty(0,T; L^2(\Omega_0^\varepsilon (0)))\cap L^2(0,T;H^1(\Omega_0^\varepsilon (0)))\right)^N\times H^1(0,T;L^2(\Omega_0^\varepsilon (0)).$$
The calculation details finally ensuring the local-in-time existence of weak solutions to the FBP are tedious. We only indicate here that the construction of the fixed-point mapping is as follows: Start off with a suitable $r^\varepsilon$ (i.e. $r^\varepsilon(t)>0$ for $t\in(0,T)$ and  $r^\varepsilon\in W^{1,\infty}_+(0,T)^M$) and solve problem (\ref{ai})--(\ref{af}). Insert the obtained solution in the equation obtained by fixing the moving boundary  in (\ref{eta2}). Refer to its solution as $s^\epsilon$.  Conclude the argument by proving for pairs of weak solutions an estimate of type $\int_0^T|s^\epsilon_1-s^\epsilon_2|^\nu\leq \lambda \int_0^T|r^\epsilon_1-r^\epsilon_2|^\nu$, where $\lambda\in (0,1)$ and $\nu\in [2,4]$. Such an estimate is a typical feature of FBPs for one-dimensional parabolic systems with kinetic  conditions at the free boundary  and can be extended to higher dimensions only if the symmetry of the evolving domains allow reduction to 1D. It shows that an improved integrability of the free boundary is available and it closes the fixed-point argument. The proof of such estimate  relies on fine energy estimates as well as on trace interpolation inequalities  combined with integral identities specific to the FBP; see, for instance, the proof of Theorem 3.4 from \cite{MB} for handling a specific 1D case and \cite{Seidman} for the 2D case.     Using Stampacchia's method in our context, suitable restrictions on data and parameters can be identified so that the obtained weak solutions stay a.e. positive.

\end{document}